\documentclass[12pt]{amsart}

\usepackage{amsaddr}
\usepackage[ansinew]{inputenc}
\usepackage[a4paper,dvips]{geometry}
\usepackage{latexsym}
\usepackage{amsfonts}
\usepackage{amsmath,amssymb}
\usepackage{bbm}
\usepackage[dvipsnames]{xcolor}
\usepackage{tikz}
\usetikzlibrary{decorations.pathreplacing,angles,quotes,calligraphy}
\usepackage{tikz-3dplot}
\usepackage{caption}
\usepackage{subfigure}
\usepackage{subcaption}
\usepackage{hyperref}

\usepackage{fixme}
\usepackage{graphicx}
\allowdisplaybreaks

\usepackage{amsmath,amssymb,amsfonts,amsthm,mathrsfs,paralist,amsrefs,fdsymbol}
\usepackage{mathtools,enumitem}
\usepackage{svg}
\usepackage{float}
\usepackage{soul}
\usepackage{tikz}
\usetikzlibrary{arrows}
\usepackage{setspace}
\usepackage{authblk}
\usepackage{dsfont}
\usepackage{tcolorbox}
\usepackage{caption,subcaption}
\usepackage{hyperref}
\theoremstyle{plain}
\newtheorem{thm}{Theorem}[section]
\newtheorem{prop}[thm]{Proposition}
\newtheorem{lemma}[thm]{Lemma}

\theoremstyle{definition}
\newtheorem{defi}[thm]{Definition}
\newtheorem{notation}[thm]{Notation}

\newtheorem{conjecture}[thm]{Conjecture}
\newtheorem{algorithm}[thm]{Algorithm}

\theoremstyle{remark}
\newtheorem{remark}[thm]{Remark}

\newtheorem{example}[thm]{Example}

\newcommand{\N}{\mathbb {N}}
\newcommand{\Z}{\mathbb {Z}}

\definecolor{Blue}{rgb}{0,0,0.7}
\newcommand{\bl}[1]{{\color{Blue}{#1}}}

\def\MSF {\mathscr{F}}

\def\BFk {\textit{\bfseries k}}

\def\iddots{\mathinner{\mkern1mu\raise1pt
    \hbox{.}\mkern2mu\raise4pt\hbox{.}\mkern2mu
        \raise7pt\vbox{\kern7pt\hbox{.}}\mkern1mu}}

\DeclareMathOperator{\supp}{supp}

\def\fooaux#1#2{%
  \mkern2mu
  \setbox0=\hbox{\mathsurround=0pt$#1\overline{\mkern-2mu #2 \mkern2mu}$}%
  \setbox1=\hbox to \wd0{\hss$#1\plus$\hss}%
  \vbox{\offinterlineskip\copy1\vskip-.6625\ht1\box0}%
  \mkern-2mu
}

\def\christiane#1{\fbox {\footnote {\ }}\ \footnotetext { From Christiane: #1}}
\def\nathan#1{\fbox {\footnote {\ }}\ \footnotetext { From Nathan: #1}}
\def\hjaspar#1{}
\def\hchristiane#1{}
\def\hnathan#1{}

\makeatletter
\g@addto@macro{\endabstract}{\@setabstract}
\newcommand{\authorfootnotes}{\renewcommand\thefootnote{\@fnsymbol\c@footnote}}%
\makeatother

\makeatletter
\def\blfootnote{\gdef\@thefnmark{}\@footnotetext}
\makeatother

\title{Golden Ratio Nets and Sequences}

\begin{document}

\maketitle

\begin{center}
  \normalsize
  \authorfootnotes
  Nathan Kirk, Christiane Lemieux\footnote{Corresponding Author} and
  Jaspar Wiart \par \bigskip

  \textit{Department of Statistics and Actuarial Science \\ University of Waterloo, Canada} \par \bigskip

  \today
\end{center}

\begin{abstract}
In this paper we introduce and study nets and sequences constructed in an irrational base, focusing on the case of a base given by the golden ratio $\varphi$. We provide a complete framework to study equidistribution properties of nets in base $\varphi$, which among other things requires the introduction of a new concept of prime elementary intervals which differ from the standard definition used for integer bases. We define the one-dimensional van der Corput sequence in base $\varphi$ and two-dimensional Hammersley point sets in base $\varphi$ and we prove some properties for $(0,1)-$sequences and $(0,m,2)-$nets in base $\varphi$ respectively. We also include numerical studies of the discrepancy of point sets and sequences in base $\varphi$ showing an improvement in distribution properties over traditional integer based Hammersley constructions. As motivation for future research, we show how the equidistribution notions that are introduced for base $\varphi$ can be generalized to other irrational bases.
\end{abstract}

\blfootnote{\textit{Contact}: {\tt{n2kirk@uwaterloo.ca, clemieux@uwaterloo.ca, jaspar.wiart@gmail.com}} }

\section{Introduction}
\label{sec:intro}

\subsection{Background \& Context}
Monte Carlo methods are used in a wide range of applications such as numerical integration and numerical approximation to solutions of differential equations, computer graphics and various aspects of mathematical finance; see \cites{DICKPILL2010, LEMIEUX2009, LEOPILL2014, NIE1992}. In the context of numerical integration, the basic idea of employing \textit{Monte Carlo integration} to approximate an integral of a function $f$ defined on $[0,1)^s$ is to collect $N$ sampling points $P_N = \{x_0, \ldots, x_{N-1}\}$ randomly in $[0,1)^s$ and determine the arithmetic mean of the $N$ values obtained as the function evaluations, i.e. $$\int_{[0,1)^s} f(x) dx \approx \frac{1}{N} \sum_{i=0}^{N-1} f(x_i).$$ The error is customarily measured via the Koksma-Hlawka inequality \cites{KOKSMA1942, HLAWKA1961} which informs us that the approximation error is at most the product between the variation of the function and the discrepancy of the $N$ sampling points. Therefore, in practise to improve the error convergence rate with $N$, one chooses the sampling points deterministically. Formally called \textit{quasi-Monte Carlo integration}, rather than choosing purely random sampling points one chooses a sequence of evaluation points which possesses a low-discrepancy and, to our good fortune, there are many good constructions known; we refer to the excellent \cite{DRMOTATICHY1997}, however by now there exists a wealth of literature on low-discrepancy point sets and sequences.

One such popular construction of low-discrepancy point sets are so-called \textit{digital nets and sequences} in an integer base $b$ which are very well studied and understood. Ever since Sobol' introduced the first construction in base $2$, in the seminal paper \cite{rSOB67a}, generalizations have been made to higher dimensions and can now be constructed for any integer base $b$. In this paper, we move away from the standard construction using an integer base, and instead work with an irrational base. While our work could be generalized to any irrational given by the largest root of $x^2-px-q$, where $1 \le q \le p$, here we focus on the case where the base is the golden ratio $\varphi = (1+\sqrt{5})/2$, i.e., the largest root of $x^2-x-1$ (when $p = q = 1$) and undertake an extensive case study on the construction of digital nets and sequence in base $\varphi$. In Section \ref{sec:otherbases} of this paper, we discuss digital nets and sequences constructed with other irrational bases.

An important aspect of the distribution properties of digital nets and sequences, and an aspect largely responsible for their initial introduction, is their intrinsic link to partitions of the unit hypercube $[0,1)^s$ into congruent hyper-rectangles which were in turn used to determine whether the construction at hand had the same number of points in each of the partitioning sets or not. Point sets which satisfy this property are called equidistributed, and this will be introduced more formally in the main body of the paper. For golden ratio based digital point sets and due to reasons involving the non-congruence of the partitioning hyper-rectangles in the golden ratio base, it is necessary to introduce a new form of equidistribution which is in some sense more granular. This new notion of \textit{strong} equidistribution is the starting point and motivation for the intuition behind point sets and sequences constructed from irrational bases possessing better distribution properties than those constructed from an integer number base.


\subsection{Previous Work}
Among the work contained in Sections \ref{sec:guidingconstructions} and \ref{sec:constructionspart2} of this text, we define and study the van der Corput sequence in base $\varphi$. It is necessary to note that studies of the one-dimensional sequence in base $\varphi$ have already been undertaken in several other works. Namely there exists two papers by S. Ninomiya \cites{Nin98, NINOMIYA1998}, the first of which defines one-dimensional sequences via $\beta-$adic transformations for $\beta \in \mathbb{R}, \beta > 1$ and using this method, the van der Corput sequence in base $\varphi$ from our Section \ref{sec:guidingconstructions} is reproduced for $\beta = \varphi$; see Example 4.1 in \cite{NINOMIYA1998}. Contained in \cite{Nin98} are numerical experiments detailing the performance of the star discrepancy for $\beta-$adic sequences for several choices of $\beta$, while the second of the Ninomiya works \cite{NINOMIYA1998} proves the exact order of the star discrepancy of the one-dimensional $\beta-$adic sequences showing that the star discrepancy of these sequences reaches the optimal order of $\log N/N$ under certain conditions on the choice of $\beta$. 

The van der Corput sequence in base $\varphi$ can also be reached via a splitting procedure of the unit interval as introduced by S. Kakutani in \cite{KAKUTANI1976}. More formally, Kakutani considers the notion of uniformly distributed sequences of nested partitions of the unit interval $[0,1)$, where at each step one refines the largest intervals in proportion $\alpha/(1-\alpha)$. Taking $\alpha = \varphi^{-1}$ and defining the points of the sequence as the end-points of the partitioning sets while imposing an ordering, the van der Corput sequence in base $\varphi$ is once again reproduced. Generalisations of Kakutani's splitting procedure have been carried out in subsequent works, and the interested reader is referred to \cite{CARBONEVOLCIC2007, VOLCIC2011}. Additionally, a specific class of refinements of partitions exists similar to those that were introduced by Kakutani. First defined by I. Carbone in \cite{CARBONE2012}, $LS-$sequences of partitions and points with respect to a real $\beta \in [0,1)$ give an explicit algorithm which orders the points determining the $LS-$sequences of partitions to obtain a uniformly distributed sequence of points, called an $LS-$sequence of points. In this context, the parameters $L$ and $S$ reference the coefficients of the order two polynomial equation $L\beta^2 + S\beta = 1$ from which the irrational base is chosen as the largest root of the equation. Clearly to obtain $\beta = \varphi$, we choose $L=1, S=-1$. Some extensions of $LS-$sequences of points to the unit square have been studied in \cite{CARBONEIACOVOLCIC2012} and a nice summary and comparison between $LS-$sequences and $\beta-$adic van der Corput sequences is given by I. Carbone in \cite{CARBONE2016}.

Most recently, and as a last comment on the existing literature, the traditional van der Corput sequence and the so-called \textit{golden ratio sequences} are implemented in a numerical comparison study regarding integration approximation contained in \cite{NIED2014}. 

\subsection{Paper Overview}
The overall contribution of this text is to formalize the study of digital nets and sequences in an irrational base with most attention given to the golden ratio base. We give a complete analogous mathematical framework to that of integer based digital nets and sequences which allows us to analyze the equidistribution properties of golden ratio based nets and sequences. Specifically, Section \ref{sec:definitionsandnotation} is used to give further context, notational conventions and several elementary facts regarding the golden ratio which will be used throughout the text. Section \ref{sec:guidingconstructions} contains definitions of the van der Corput sequence and Hammersley point set in base $\varphi$ introducing two important constructions which will provide motivation and direction for our study moving forward. The framework surrounding equidistribution properties of golden ratio based nets and sequences is developed in Sections \ref{sec:elemint} and \ref{sec:equid} with several results regarding the equidistribution of the van der Corput sequence, Hammersley point sets and more general $(t,m,s)-$nets and $(t,s)-$sequences in base $\varphi$ are proven in Section \ref{sec:constructionspart2}. Also contained in Section \ref{sec:constructionspart2}, we discuss some fundamental obstacles with randomizing point sets in base $\varphi$ in the form of scrambling while also providing some recommendations for future work in this area. of Section \ref{sec:numericalresults} contains numerical calculations of the star discrepancy of the Hammersley point set in base $\varphi$ showing a significant improvement in distribution properties when comparing to classical integer based Hammersley constructions. To finish, Section \ref{sec:otherbases} gives further generalized notions of equidistribution of point sets and sequences constructed from an irrational base $\gamma \in \mathbb{R}$ acting as the largest root from the polynomial $$x^2-px-q$$ for $1\leq p\leq q$. There are two appendices to this text which gives further supporting graphics of irrational based point sets and tables of numerical results on the star and $L_2-$discrepancy for completeness.

\section{Basic Notation and Definitions}\label{sec:definitionsandnotation}

To begin, we present elementary facts and relations of the golden ratio and give notational conventions that will be used throughout. For the majority of this paper, we study the point sets and sequences generated by utilizing the famous golden ratio $\varphi$ as a base. The golden ratio is a particularly nice irrational number to use as a base due to the fact that each $x \in \mathbb{R}$ has a \textit{unique}, possibly infinite, base $\varphi$ expansion which we will denote as
\[
x=\sum_{j=m-1}^{-\infty} d_j\varphi^j=(d_{m-1}\dots d_1d_0.d_{-1 }d_{-2}\dots)_\varphi
\]
where for all valid $j$ we have $d_j\in\{0,1\}$ and $d_jd_{j-1}=0$ (i.e., there are no consecutive ones). We can guarantee this \textit{reduced form} representation since from the minimal polynomial, $\varphi$ satisfies $\varphi^2 = \varphi + 1$, i.e. we can replace the string $(011)_\varphi$ by $(100)_{\varphi}$ in base $\varphi$ expansions. We make a few initial observations:
\begin{itemize}
    \item Note that, e.g.,  $5 = \varphi^3 + \varphi^{-1} + \varphi^{-4} = (1000.1001)_\varphi$ and $8 = \varphi^4 + \varphi^0 + \varphi^{-4} = (10001.0001)_\varphi$ have finite representations in base $\varphi$. But $\pi = (100.010010101001...)_\varphi$ has an infinite base $\varphi$ representation (sequence \href{https://oeis.org/A102243}{A102243} in the \href{https://en.wikipedia.org/wiki/On-Line_Encyclopedia_of_Integer_Sequences}{OEIS}).
    \item The set of numbers with finite decimal expansion in base $\varphi$ is the ring $\Z[\varphi]$. 
    \item Finally, examples of adding and subtracting numbers in base $\varphi$ are
    \[
    2=(1)_\varphi+(1)_\varphi = (1)_\varphi+(0.11)_\varphi=(1.11)_\varphi=(10.01)_\varphi.
    \]
    and
    \[
    \varphi - 1 = (10)_\varphi-(1)_\varphi = (1.1)_\varphi-(1)\varphi = (0.1)_\varphi = 1/\varphi.
    \]
\end{itemize}

The golden ratio possesses an intrinsic link with the famous Fibonacci sequence of integers which will play a major role throughout this paper.

\begin{defi}
The \emph{Fibonacci sequence} is defined by $F_0=0$, $F_1=1$, and $F_k=F_{k-1}+F_{k-2}$ for all $k\geq 2$. It will be convenient for us to use the notation $F^m=F_{m+2}$ where $m \ge -2$.
\end{defi}

In addition, the Fibonacci numbers serendipitously work as a base for $n\in\N_0$, which we will refer to as base $\MSF$. Specifically, non-negative integers $F^{m-1}\leq n<F^m$ can be written in base $\MSF$ as 
\[
n=(d_{m-1}\dots d_1 d_0)_\MSF=\sum_{j=0}^{m-1}d_jF^j,
\]
where there does not exist $0<i<m$ such that $d_{i}=d_{i-1}=1$, i.e., once again there are no consecutive ones. This can be shown via a simple induction argument. Indeed, assume that all integers $0\leq k<F^{m}$ can be written in this way without consecutive ones and choose $F^{m}\leq n <F^{m+1}$. Then $0\leq n-F^{m} \leq F^{m-2}$ and $n=(10d_{m-2}\dots d_1 d_0)_\MSF$, where $n-F^m=(d_{m-2}\dots d_1 d_0)_\MSF$. In this paper we will primarily be concerned with whole numbers in base $\varphi$ for which the following notation will be helpful.
\begin{notation}
    For $n=(d_{m-1}\dots d_1d_0)_\MSF\in\N_0$, we use the notation $\overline n$ to denote $(d_{m-1}\dots d_1d_0)_\varphi$, the $n^{th}$ whole number in base $\varphi$. The zeroth digit of either $n$ or $\overline n$ is denoted by $|n|=|\overline n|=d_0$. For $k\in\Z$ we also define the operators
    \[
    n\odot\varphi^k=\sum_{j=0}^{m-1}d_jF^{j+k}
    \quad\text{and}\quad
    \overline n\oplus k=\overline{n+k}
    \]
    and note that in this paper we only use those $k$ such that $n\odot\varphi^k\in\N_0$ or $n+k\in\N_0$.
\end{notation}

One of the features of base $\varphi$ is that the difference between two consecutive whole numbers in base $\varphi$ could be either $1$ or $\varphi^{-1}$. This solely depends on the zeroth digit of $\overline{n}$, $d_0$. 
\begin{lemma}\label{lem: consecutive difference}
    The difference between a whole number $\overline{n} = (d_{m-1}\dots d_2 d_1 d_0)_\varphi$ in base $\varphi$ and the next is 
    \[
    \overline{n+1}-\overline n = \varphi^{-|\overline{n}|}.
    \]
\end{lemma}
\begin{proof}
To verify this relation, simply evaluate the difference $\overline{n+1}-\overline n$ in each of the following exhaustive cases
\begin{enumerate}
    \item [(i)] $n=(d_{m-1}\dots d_2 00)_\MSF$ so that $n+1=(d_{m-1}\dots d_201)_\MSF$,
    \item[(ii)] $n=(d_{m-1}\dots d_2 01)_\MSF$ so that $n+1=(d_{m-1}\dots d_210)_\MSF$, and 
    \item[(iii)] $n=(d_{m-1}\dots d_3010)_\MSF$ so that $n+1=(d_{m-1}\dots d_3100)_\MSF$
\end{enumerate}
The result follows by applying the relation $\varphi^2 = \varphi+1$ when necessary.
Note that $n+1$ might not be written in its reduced form in the above so that there might be consecutive ones. However, this does not affect our conclusion or the statement of the lemma.
\end{proof}

To conclude this section, we recall several identities involving the golden ratio and Fibonacci numbers which will, at times, be useful moving forward:
\begin{align}
    F_{m}&=F^{m-2} = \frac{\varphi^m-\psi^m}{\sqrt{5}}, m \ge 0,  \text{ where } \psi=1-\varphi=-\varphi^{-1}\label{eq: fib as phi}\\
    \varphi^{m}&=F_{m}\varphi+F_{m-1} = F^{m-2}\varphi+F^{m-3}, m \ge 1, \text{ and}\label{eq: varphi as fib}\\
    \psi^m&=F_m\psi +F_{m-1}=F^{m-2}\varphi+F^{m-3}, m \ge 1.
\end{align}

\section{Guiding examples for constructions}\label{sec:guidingconstructions}

Before going any further, it is useful to describe the basic constructions that will be used as building blocks in our framework. We consider prototypical constructions of low discrepancy point sets in one and two dimensions. In dimension one, the archetypal sequence example in base $\varphi$ is the \emph{van der Corput sequence} which we now define in the usual way as when considering an integer base.

\begin{defi}
    If $n=(d_{m-1}\dots d_1 d_0)_\MSF\in\N_0$, then the $n^{th}$ point of the \textit{van der Corput sequence in base $\varphi$} is given by $g_n := (.d_0 d_1\dots d_{m-1})_\varphi$. 
\end{defi}

\begin{figure}[h!]
    \centering
    \begin{tikzpicture}

        \draw[thick, black] (0,0) -- (10,0);

        \draw[thin,black] (0,0.2) -- (0,-0.2);
        \draw[thin,black] (10,0.2) -- (10,-0.2);
        
        \filldraw[black] (0,0) circle (3pt);
        \filldraw[black] (6.180,0) circle (3pt);
        \filldraw[black] (3.819,0) circle (3pt);
        \filldraw[black] (2.360,0) circle (3pt);
        \filldraw[black] (8.541,0) circle (3pt);
        \filldraw[black] (1.458,0) circle (3pt);
        \filldraw[black] (7.639,0) circle (3pt);
        \filldraw[black] (5.278,0) circle (3pt);
        \filldraw[black] (0.901,0) circle (3pt);
        \filldraw[black] (7.082,0) circle (3pt);
        \filldraw[black] (4.721,0) circle (3pt);
        \filldraw[black] (3.262,0) circle (3pt);
        \filldraw[black] (9.442,0) circle (3pt);

        \node at (0,-0.5) {$0$};
        \node at (10,-0.5) {$1$};
        \node at (6.180, -0.5) {$\varphi^{-1}$};
        \node at (3.819,-0.5) {$\varphi^{-2}$};
        \node at (2.36,-0.5) {$\varphi^{-3}$};
        \node at (1.458,-0.5) {$\varphi^{-4}$};

    \end{tikzpicture}
    \caption{The first 13 terms in the van der Corput sequence in base $\varphi$}
    \label{fig:vdcbasevarphi - first 13 terms}
\end{figure}
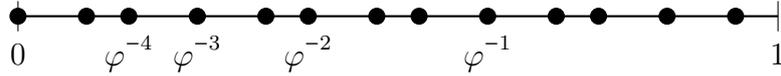

The first few points of this sequence are $g_0 = (.0)_\varphi = 0$, $g_1 = (.10)_\varphi = \varphi^{-1}$, $g_2 = (.010)_\varphi = \varphi^{-2}$, $g_3 = (.0010)_\varphi = \varphi^{-3}$, $g_4 = (.101)_\varphi = \varphi^{-1}+\varphi^{-3}$, $g_5 = (.0001)_\varphi = \varphi^{-4}$, $g_6 = (.1001)_\varphi = \varphi^{-1}+\varphi^{-4}$, $g_7 = (.0101)_\varphi = \varphi^{-2}+\varphi^{-4}$, etc. For illustration purposes, Figure \ref{fig:vdcbasevarphi - first 13 terms} gives the first $F^5 = 13$ terms of the van der Corput sequence in base $\varphi$. It is also worthwhile to refer to Figure \ref{fig:elem m-int} which shows the embedding of the terms of the van der Corput sequence in base $\varphi$ for each $F^m$ number of terms for $0 \leq m \leq 5$.

In two dimensions, we study the \textit{Hammersley point set}. One way to describe the Hammersley construction in base $\varphi$ (and also in the integer base) is to form the $i^{th}$ two-dimensional point by keeping the respective term of the van der Corput sequence as the first coordinate, and adding a second coordinate obtained by reversing the order of the digits of the number defining the first coordinate. Alternatively, the second coordinate in the $i^{th}$ point can also be seen to be the ``$i^{th}$ whole number" (in a given base) divided by the total number of points in the set. We will refer back to this construction algorithm for the Hammersley set in Section \ref{sec:otherbases}. Therefore, in our case of the golden ratio and due to the fact that for a given $m \geq 0$ there are $F^m$ $m-$digit whole numbers in base $\MSF$, we fix the number of points in advance to be $F^m$.

\begin{defi}
Fix $m \in \mathbb{N}$. The \textit{Hammersley point set in base $\varphi$} with $F^m$ points is defined as
\[
H_m = \big\{ \left( (.d_0d_1\dots d_{m-1})_\varphi,(.d_{m-1}\dots d_1d_0)_\varphi \right): 0\leq (d_{m-1}\dots d_1d_0)_\MSF<F^m\big\}.
\]
\end{defi}

\begin{figure}[t!]
     \centering
     \includegraphics[width=14cm]{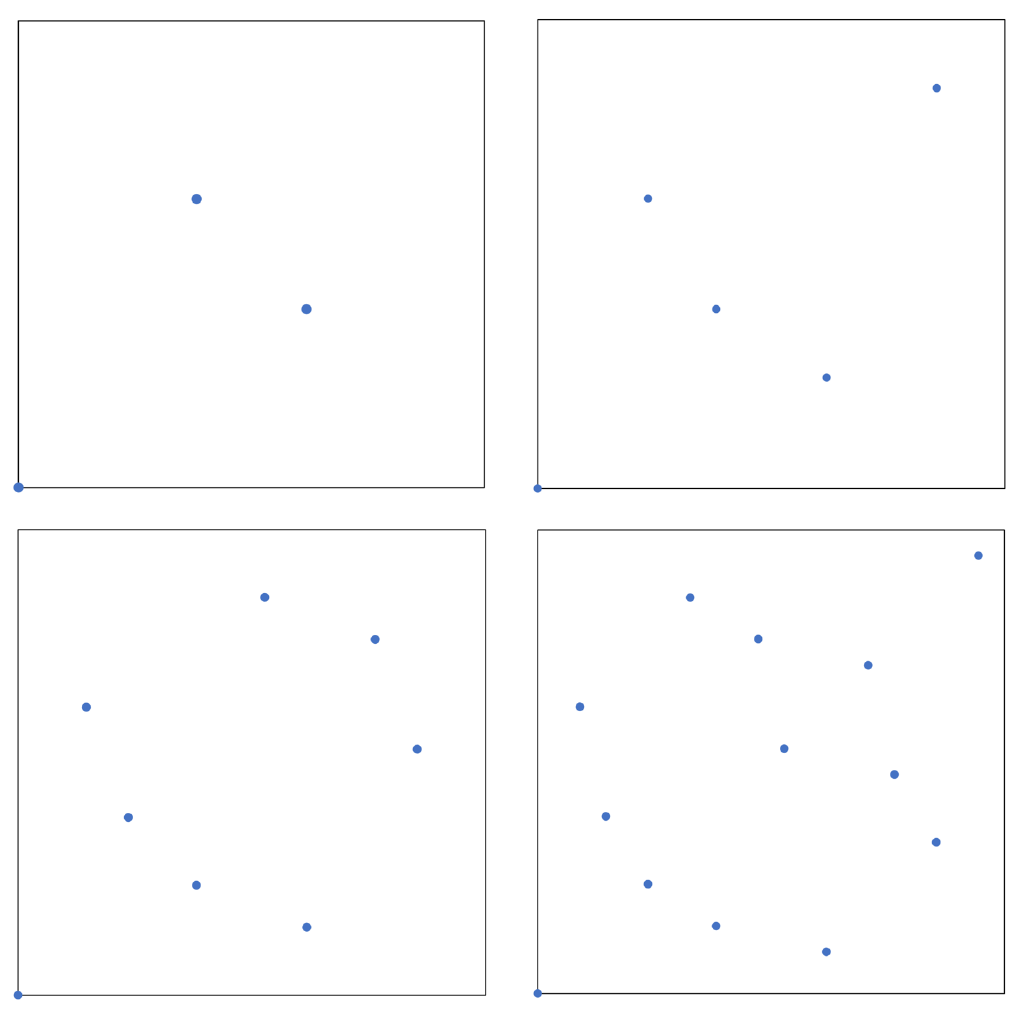}
     \caption{The Hammersley point set in base $\varphi$ for $m = 2, 3, 4$ and $5$.}\label{fig:Hammersleypointsets1-4}
\end{figure}

For example, when $m=3$ we have that
\[
H_3 = \left\{(0,0),((0.100)_\varphi,(0.001)_\varphi),((0.010)_\varphi,(0.010)_\varphi),((0.001)_\varphi,(0.100)_\varphi),((0.101)_\varphi,(0.101)_\varphi)\right\}
\]
and for the reader's benefit, $H_m$ is given for $m=2, 3, 4$ and $5$ as Figure \ref{fig:Hammersleypointsets1-4} and $H_{16}$ is given as Figure \ref{fig: hammersley 16} in Appendix \ref{app:additional figures 1}.

\section{Elementary intervals}
\label{sec:elemint}

In the literature, the concept of an elementary interval can be traced back to the seminal paper \cite{rSOB67a} from 1967. The author, I. M. Sobol', at that point only working in base $2$, referred to the modern elementary intervals in base $2$ as \textit{binary segments} in one dimension and \textit{binary parallelepipeds} in dimension greater than one. More than a decade later, the notion was generalized by H. Faure to bases $b\geq 2$ in \cite{FAURE1981}, however it wasn't until a follow-up publication \cite{FAURE1982} in which these intervals were called \textit{``pav\'{e}s \'{e}l\'{e}mentaires"}. All of these authors were inspired by studying the basic distribution properties of the classical van der Corput sequence and its link to partitions of the unit interval $[0,1)$. For the benefit of the reader, we briefly recall the definitions of elementary intervals, equidistribution and nets and sequences when working with integer bases. For $b \geq 2$, an \emph{elementary $k$-interval in base $b$} is a subset of $[0,1)$ of the form
\[
\left[\frac{a}{b^k},\frac{a+1}{b^k}\right)
\]
where $0\leq a<b^k$ and $k \in \mathbb{N}_0$. In $s-$dimensions, an \emph{elementary $(k_1,\dots,k_s)-$interval in base $b$} is the product of elementary $k_j-$intervals in base $b$, i.e. it is a set of the form
\[
\prod_{j=1}^s\left[\frac{a_j}{b^{k_j}},\frac{a_j+1}{b^{k_j}}\right),
\]
where $0\leq a_j<b^{k_j}$ for each $1\leq j \leq s$. For $m\geq 0$, a point set $P_m$ contained in $[0,1)^s$ of size $b^m$ is called \emph{$(k_1,\dots,k_s)-$equidistributed in base $b$} if every elementary $(k_1,\dots,k_s)-$interval in base $b$ contains exactly $b^{m-k_1-\dots-k_s}$ points from $P_m$. The point set $P_m$ containing $b^m$ points is called a \emph{$(t,m,s)-$net in base $b$} if it is $(k_1,\dots,k_s)-$equidistributed in base $b$ for all $k_1+\dots +k_s \leq m-t$. Finally, a sequence $\{x_0,x_1,\dots\}\subset[0,1)^s$ is called a \emph{$(t,s)-$sequence in base $b$} if for every $m>t$ and $k>0$, the set $\{x_{kb^m},\dots,x_{kb^m+b^m-1}\}$ forms a $(t,m,s)-$net.

\begin{figure}[t!]
\begin{center}
\includegraphics[width=\textwidth]{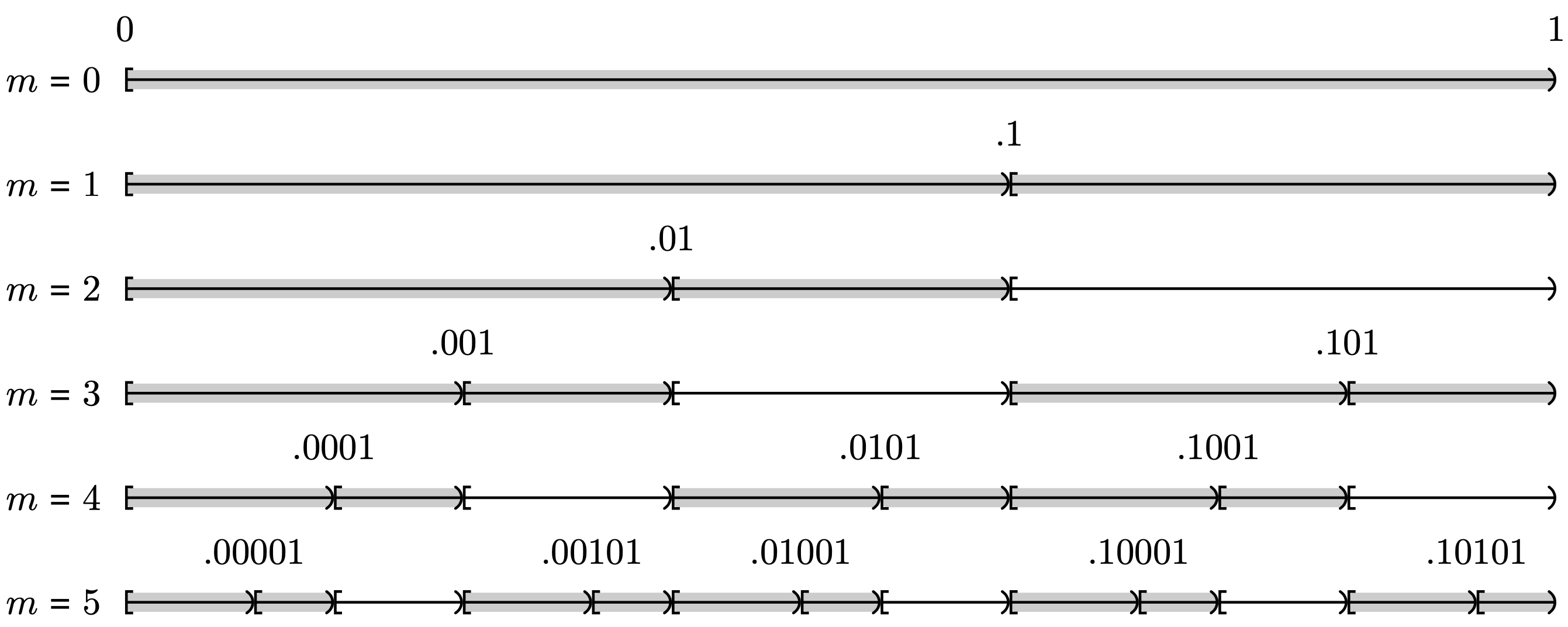}
\caption{The first six $m$-partitions $\mathcal{P}_m$, with the prime elementary intervals filled in gray.}
\label{fig:elem m-int}
\end{center}
\end{figure}

The goal for this section is to develop the definition of an elementary interval in base $\varphi$ from which everything else regarding the distribution notions will automatically follow. To begin forming the analogous definitions in base $\varphi$, it is clear that an elementary interval in base $\varphi$ should be a subset of $[0,1)$ of the form 
\[
I=\left[\frac{\overline{a}}{\varphi^m},\frac{\overline{a+1}}{\varphi^m}\right),
\]
\hchristiane{Maybe we should first five def 4.1 and then explain. Think about this.}
where $0\leq a<F^m$ and $m \in \mathbb{N}_0$. It is clear that for a fixed $m$, the set of all such $I$ forms a partition of $[0,1)$ which we call the \emph{$m-$partition} and denote this object by $\mathcal{P}_m$. The first six such partitions are given in Figure \ref{fig:elem m-int}, and from this illustration we can observe several extremely important properties which must be taken into consideration for the remainder of this text. 

Unlike when forming elementary intervals for an integer base, for $m\geq1$ the elementary intervals in base $\varphi$ in $\mathcal{P}_m$ have two different lengths of $\varphi^{-m}$ and $\varphi^{-m-1}$. We refer back to Lemma \ref{lem: consecutive difference} where we saw that $\overline{a+1}-\overline{a}=\varphi^{-|a|}$ can be either $1$ or $\varphi^{-1}$, i.e., the size of an elementary interval $I$ in base $\varphi$ depends upon the numerator $0 \leq a < F^m$. Second, we see that $\mathcal{P}_{m+1}$ refines $\mathcal{P}_m$ with the elementary intervals of length $\varphi^{-m}$ being split into two subintervals of lengths $\varphi^{-m-1}$ and $\varphi^{-m-2}$, while the elementary intervals of length $\varphi^{-m-1}$ are not refined and remain unchanged in $\mathcal{P}_{m+1}$. In other words, the largest intervals in a given partition are refined in the next partition, and the smallest remain the same. Third, some of the elementary intervals in $\mathcal{P}_m$ appear in $\mathcal{P}_{m-1}$ while others appear in $\mathcal{P}_{m+1}$. This last observation is key and means that if we were to treat all elementary intervals in $\mathcal{P}_m$ equally, then there would be some elementary $m-$intervals that are also elementary $(m-1)-$intervals and others that would also be elementary $(m+1)-$intervals. So the question arises, how do we categorize and differentiate between elementary intervals of the same size which appear in more than one $m-$partition of the unit interval? With this in mind, the following definition is natural.

\begin{defi}
    An elementary $m-$interval $I$ in base $\varphi$ is a subinterval of $[0,1)$ of the form $$I = \left[\frac{\overline a}{\varphi^m},\frac{\overline{a+1}}{\varphi^m}\right)$$ for some $0 \leq a <F^m$ and $m \in \mathbb{N}_0$.
\end{defi}

When convenient to do so, we shall simply refer to an elementary $m-$interval as an elementary interval. We usually do so when the degree of refinement of the elementary interval is not important for the discussion. To address the problem observed in the last paragraph, that elementary intervals of the same size appear in multiple levels of partition, we now distinguish between intervals by introducing the new notion of a \textit{prime} elementary interval.

\begin{defi}
An elementary interval $I$ in base $\varphi$ is called a \emph{prime elementary $k-$interval in base $\varphi$} if $k$ is the least positive integer such that there exists $0\leq a<F^k$ such that
\[
I=\Bigg[\frac{\overline a}{\varphi^k},\frac{\overline{a+1}}{\varphi^k}\Bigg).
\]
\end{defi}

\begin{remark}
    To further help the reader distinguish between elementary intervals and prime elementary intervals, we use the exponent notation $m$ for a general elementary interval contained in the $m-$partition ${\mathcal P}_m$, $$\left[\frac{\overline a}{\varphi^m},\frac{\overline{a+1}}{\varphi^m }\right)$$ for some $0 \leq a \leq F^m$, and we use the notation $k$ for the base exponent when the elementary interval is a prime elementary interval.
\end{remark}

We wish to characterize all prime elementary intervals and for this purpose, the following notation will be useful.

\begin{notation}\label{not:elementary intervals}
Let $I$ be an elementary $m-$interval with
\[
I=\Bigg[\frac{\overline b}{\varphi^m},\frac{\overline{b+1}}{\varphi^m}\Bigg),
\]
where $0\leq b=(d_{m-1}\dots d_1d_0)_\MSF<F^m$. Set $a=(d_{m-1}\dots d_200)_\MSF$ and choose $i\in\{0,1,2\}$ such that $b=a+i$. We use the notation $I_i(a;m)$ to denote $I$, i.e.\
\begin{equation}
\label{eq:Iiam}
I=I_i(a;m):=\Bigg[\frac{\overline {a+i}}{\varphi^m},\frac{\overline{a+i+1}}{\varphi^m}\Bigg).
\end{equation}
\end{notation}
\hchristiane{I like the notation, everything is clear. I wonder if it would help the reader to say $i$ is the type, $a$ is some kind of anchor, and $m$ refers to the level of the partition}

We emphasise the following important aspects of this notation which is meant to keep track of the last two digits of $b$ in the above notation. Indeed if $d_1=d_0=0$ then $i=0$, if $d_1=0$ and $d_0=1$ then $i=1$, and if $d_1=1$ and $d_0=0$ then $i=2$. Hence, in Lemma \ref{lem:1d elem int characterization} we characterize the elementary intervals in $\mathcal{P}_m$ in terms of $i$ and we can think of $i$ as the \textit{type} of the elementary interval $I\subset \mathcal P_m$ with respect to the $m-$partition. Importantly, if $I$ is contained in $\mathcal{P}_{m-1}$ or $\mathcal{P}_{m+1}$, then it will have a different type with respect to those partitions.  As a second note, when $a=(d_{m-1}\dots d_3000)_\MSF$ then $I_0(a;m)$, $I_1(a;m)$ and $I_2(a;m)$ all exist, however when $a=(d_{m-1}\dots d_40100)_\MSF$ only $I_0(a;m)$ and $I_1(a;m)$ exist.

The following lemma states a useful property about the length of the interval $I_i(a;m)$ depending on the value of $i$.
\begin{lemma}\label{lem:length of interval}

    The length of $I_i(a;m)$ is $\varphi^{-m}$ when $i=0$ or $i=2$ and is $\varphi^{-m-1}$ when $i=1$.
\end{lemma}

\begin{proof}
This follows directly from Lemma \ref{lem: consecutive difference}.
\end{proof}

We now prove the second and third observations we made about $\mathcal{P}_m$ in the above paragraph as well as characterizing prime elementary intervals in base $\varphi$.
\begin{lemma}\label{lem:1d elem int characterization}
Let $m\in \N_0$ and $0\leq a<F^m$. 
\begin{enumerate} 
    \item[0.]  Elementary intervals $I_0(a;m)$ appear only in $\mathcal{P}_m$ and are refined as 
    \[
    I_0(a;m)=I_0(a\odot\varphi;m+1)\cup I_1(a\odot\varphi;m+1)
    \]
    in $\mathcal{P}_{m+1}$. For $m>0$, there are $F^{m-2}$ elementary intervals of this form.
    \item[1.] Elementary intervals $I_1(a;m)$ appear in $\mathcal{P}_m$ and $\mathcal{P}_{m+1}$ only as 
    \[
    I_1(a;m)=I_2(a\odot\varphi;m+1).
    \]
    For $m>0$, there are $F^{m-2}$ elementary intervals of this form.
    \item[2.] Elementary intervals $I_2(a;m)$ appear in $\mathcal{P}_m$ and $\mathcal{P}_{m-1}$ only and are refined in $\mathcal{P}_{m+1}$ as:
    \[
    I_2(a;m)= I_1(a\odot \varphi^{-1};m-1)=I_0(a\odot\varphi+3;m+1)\cup I_1(a\odot\varphi+3;m+1).
    \]
    For $m>0$, there are $F^{m-3}$ elementary intervals of this form.
\end{enumerate}
 In particular there are $F^m$ elementary intervals in $\mathcal{P}_m$ and an elementary interval 
\[
I=\Bigg[\frac{\overline a}{\varphi^m},\frac{\overline{a+1}}{\varphi^m}\Bigg),
\]
with $0\leq a=(d_{m-1}\dots d_1d_0)_\MSF<F^m$, is a prime elementary $k-$interval if and only if $d_1=0$, i.e.\ if and only if $I=I_0(a;k)$ or $I=I_1(a;k)$.
\end{lemma}

At this point, it may be helpful for the reader to study Figure \ref{fig:elem m-int} which distinguishes the prime elementary intervals within the $m-$partition for $0 \leq m\leq 5$.

\begin{proof}[Proof of Lemma \ref{lem:1d elem int characterization}]
The main idea of this proof is to compare the lengths  of elementary intervals in consecutive partitions with the same left endpoint to determine if they are equal. We first consider elementary intervals in $\mathcal{P}_m$ of the form 
\[
I_0(a;m)=\Bigg[\frac{\overline {a}}{\varphi^m},\frac{\overline{a+1}}{\varphi^m}\Bigg),
\]
where $0\leq a=(d_{m-1}\dots d_200)_\MSF<F^m$. From Lemma \ref{lem:length of interval} we know that the length of these elementary intervals is $\varphi^{-m}$, and since $\mathcal{P}_m$ only has elementary intervals of length $\varphi^{-m}$ and $\varphi^{-m-1}$, the only other partition that the elementary interval in question can be a part of is $\mathcal{P}_{m-1}$. The elementary interval in the $\mathcal{P}_{m-1}$ whose left endpoint coincides with that of $I_0(a;m)$ is
\[
I=\Bigg[\frac{\overline{ a\odot\varphi^{-1}}}{\varphi^{m-1}},\frac{\overline{a\odot\varphi^{-1}+1}}{\varphi^{m-1}}\Bigg),
\]
but $\overline{a\odot\varphi^{-1}}$ ends in $0$, so the length of these elementary intervals is $\varphi^{-m+1}$, and therefore $I_0(a;m)$ only appears in $\mathcal{P}_{m}$. To see how $I_0(a;m)$ splits in $\mathcal{P}_{m+1}$, we simply note that 
\[
\frac{\overline{a\odot\varphi}}{\varphi^{m+1}}=\frac{\overline a}{\varphi^m}
\]
which means that $I_0(a\odot\varphi;m+1)$ has the same left endpoint as $I_0(a;m)$. This, combined with the fact that the length of $I_0(a\odot\varphi;m+1)\cup I_1(a\odot\varphi;m+1)$ is $\varphi^{-m-1}+\varphi^{-m-2}=\varphi^{-m}$ tells us that $I_0(a;m)$ splits in $\mathcal{P}_{m+1}$ into these two elementary $(m+1)-$intervals.

Next, we consider elementary intervals in $\mathcal{P}_{m}$ of the form
\[
I_1(a;m)=\Bigg[\frac{\overline {a+1}}{\varphi^m},\frac{\overline{a+2}}{\varphi^m}\Bigg),
\]
where $0\leq a=(d_{m-1}\dots d_200)_\MSF<F^m$. Elementary intervals of this kind have length $\varphi^{-m-1}$ and with similar reasoning as above we conclude that the only other partition which could contain these elementary intervals is $\mathcal{P}_{m+1}$. Since $(a+1)\odot\varphi=(d_{m-1}\dots d_{3}010)_\MSF=a\odot\varphi+2$, we see that
\[
\frac{\overline{a+1}}{\varphi^{m}}=\frac{\overline{(a+1)\odot\varphi}}{\varphi^{m+1}}=\frac{\overline{a\odot\varphi+2}}{\varphi^{m+1}}.
\]
We conclude that $I_1(a;m)=I_2(a\odot\varphi;m+1)$ because they have the same left endpoint and the same length.

The last of the elementary intervals in $\mathcal{P}_{m}$ which we need to consider are the form
\[
I_2(a;m)=\Bigg[\frac{\overline {a+2}}{\varphi^m},\frac{\overline{a+3}}{\varphi^m}\Bigg),
\]
where $0\leq a=(d_{m-1}\dots d_3000)_\MSF<F^m$. As above, these intervals have length $\varphi^{-m}$ and the only possible partition besides $\mathcal{P}_{m}$ which can contain these elementary intervals is $\mathcal{P}_{m-1}$. Since $(a+2)\odot\varphi^{-1}=(d_{m-1}\dots d_{3}01)_\MSF=a\odot\varphi^{-1}+1$, the left endpoint of $I_2(a;m)$ and $I_1(a\odot\varphi^{-1};m-1)$ coincide. We conclude that these two elementary intervals are the same because they have the same left endpoint and the same length. Now $(a+2)\odot\varphi=(d_{m-1}\dots d_30100)_\MSF=a\odot\varphi+3$ so that $I_2(a;m)$ and $I_0(a\odot\varphi+3;m+1)$ have the same left endpoint. Since the length of $I_2(a;m)$ and $I_0(a\odot\varphi+3;m+1)\cup I_1(a\odot\varphi+3;m+1)$ are the same, $I_2(a;m)$ splits in $\mathcal{P}_{m+1}$ exactly as claimed.

We count the elementary intervals of each type by induction. The base case is $\mathcal{P}_{1}$ which can easily be verified in Figure \ref{fig:elem m-int}. In that partition, there is $F^{-1}=1$ elementary interval of type $i=0$, namely $I_0(0;1)$, there is $F^{-1}=1$ elementary interval of type $i=1$, namely $I_1(0;1)$, and $F^{-2}=0$ elementary intervals of type $i=2$ as claimed. Suppose that there are $F^{m-2}$ elementary intervals of type $i=0$ and $i=1$ and $F^{m-3}$ elementary intervals of type $i=2$ in $\mathcal{P}_{m}$. As we have shown, every elementary interval of type $i=0$ and $i=2$ splits into both an elementary interval of type $i=0$ and $i=1$ in $\mathcal{P}_{m+1}$ and every elementary interval of type $i=1$ in $\mathcal{P}_{m}$ becomes an elementary interval of type $i=2$ in $\mathcal{P}_{m+1}$. Thus there are $F^{m-2}+F^{m-3}=F^{(m+1)-2}$ elementary intervals of type $i=0$ and $i=1$ in $\mathcal{P}_{m+1}$ and $F^{m-2}=F^{(m+1)-3}$ elementary intervals of type $i=2$ in $\mathcal{P}_{m+1}.$

In total, there are $F^{m-2}+F^{m-2}+F^{m-3} = F^m$ elementary intervals in $\mathcal{P}_m$.
\end{proof}

Now that we have defined and classified all prime elementary intervals, the remaining definitions follow without much difficulty. We first make one small note, which will play a larger role in Section \ref{sec:equid} concerning equidistribution. Earlier in this section, we defined the $m-$partition of $[0,1)$ and assigned the notation $\mathcal{P}_m$ for $m \in \mathbb{N}_0$. Analogously in $s$ dimensions, for a vector $\mathbf{m} = (m_1, \ldots, m_s) \in \mathbb{N}_0^s$ we define the $(m_1, \ldots, m_s)-$partition of $[0,1)^s$ in the natural way (i.e. the $j^{th}$ axis is partitioned into $\mathcal{P}_{m_j}$) and denote this partition $\mathcal{P}_{\mathbf{m}}$.

\begin{defi}
    An \textit{elementary $(m_1, \ldots, m_s)-$interval $I$ in base $\varphi$} is a subinterval of $[0,1)^s$ of the form $$I = \prod_{j=1}^s \left[ \frac{\overline{a_j}}{\varphi^{m_j}}, \frac{\overline{a_j+1}}{\varphi^{m_j}} \right)$$ where $0 \leq a_j \leq F^{m_j}$ and $m_j \in \mathbb{N}_0$ for each $j$.
\end{defi}

\begin{defi}

    A \emph{prime elementary $(k_1,\dots,k_s)-$interval in base $\varphi$} is a subset of $[0,1)^s$ of the form
    \[
    I=\prod_{j=1}^s\Bigg[\frac{\overline a_j}{\varphi^{k_j}},\frac{\overline{a_j+1}}{\varphi^{k_j}}\Bigg)
    \]
    where $0\leq a_j<F^{k_j}$ has its second-to-last digit equal to zero in its base $\MSF$ expansion, and $k_j \in \mathbb{N}_0$ for each $j$. 
\end{defi}  

\begin{figure}[t!]
     \centering
     \includegraphics[scale = .45]{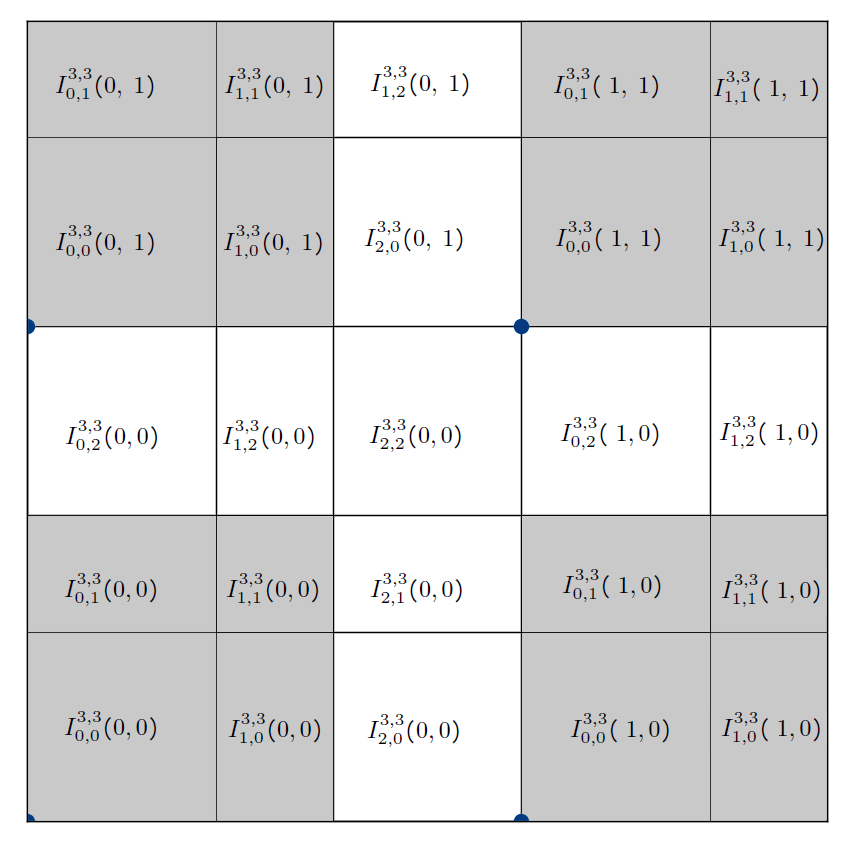}
     \caption{All elementary intervals in the $(3,3)$-partition with the prime elementary $(3,3)$-intervals shaded in gray; these come in groups of four. We have written $I_{i,j}(a,b;m_1,m_2)$ compactly as $I_{i,j}^{m_1,m_2}(a,b)$.}\label{fig: 2d elem int type}
\end{figure}

As a reference, Figure \ref{fig: 2d elem int type} gives an example of a partition $\mathcal{P}_{(3,3)}$ in two dimensions with the higher dimensional prime elementary intervals shown once again in gray. Moreover, one can visualize that in two dimensions, we have four differing sizes of elementary intervals and in general, in $s$ dimensions we will have a possible $2^s$ different sizes of elementary intervals in $\mathcal{P}_{\mathbf m}$ for $\mathbf{m} \in \mathbb{N}_0^s$ with the prime elementary intervals forming groups of $2^s$. Using Notation \ref{not:elementary intervals}, we note that prime elementary intervals are products of intervals of the form $I_i(\cdot;\cdot)$ with $i=0,1$ and not $i=2$. Formally, in two-dimensions for $\mathbf{k}= (k_1,k_2)\in\N^2_0$ we define
\[
I_{i,j}(a,b;k_1,k_2):=I_i(a;k_1)\times I_j(b;k_2)=\Bigg[\frac{\overline {a+i}}{\varphi^{k_1}},\frac{\overline{a+i+1}}{\varphi^{k_1}}\Bigg)\times\Bigg[\frac{\overline {b+j}}{\varphi^{k_2}},\frac{\overline{b+j+1}}{\varphi^{k_2}}\Bigg).
\]
As in the one-dimensional case for a given $a$ or $b$ whose last two digits are zero, the case $i=2$ or $j=2$ might not be defined. Lastly, when we do not need to be specific to which prime elementary interval we are referring, we will simply write 
\begin{equation}\label{not:arbprimeinterval}
I(a,b;k_1,k_2)
\end{equation}
to denote one of the four prime elementary intervals $I_{i,j}(a,b;k_1,k_2)$ with $i,j\in\{0,1\}$.

To finish this section we establish some notation which is important for determining the volume of a given $s-$dimensional elementary interval in base $\varphi$.

\begin{notation}    
For an elementary interval of the form $$ I = \prod_{j=1}^s \left[ \frac{\overline{a_j}}{\varphi^{m_j}}, \frac{\overline{a_j+1}}{\varphi^{m_j}} \right)$$ where $0\leq a_j<F^{m_j}$, we will use the notation 
    \begin{equation}
    \label{eq:logVolumeI}
        |I|=\sum_{j=1}^s \left( m_j+|a_j|\right)
    \end{equation}
    so that the volume of $I$ is $\varphi^{-|I|}$. Take care in noting that the action of the notation $|\cdot|$ depends upon whether the argument is an interval $I$, or a digit $a_j$ recalling that $|a_j|$ refers to the zeroth digit of $a_j$. Further note that in the one-dimensional case Definition \ref{eq:logVolumeI} is equivalent to the existing definition for the length of an elementary interval as given in Lemma \ref{lem:length of interval}.
\end{notation}
\hchristiane{We may need to be more precise with the notation used in the above definition. E.g., ok to use $d_1$ without a $j$? Ok to use $a$ rather than $b$ from definition of $I_i(a;m)$?}

\section{Equidistribution Definitions and Properties}\label{sec:equid}

Using the definition of prime elementary $(k_1, \ldots, k_s)-$intervals in base $\varphi$ from the previous section, we can now define  $(k_1,\dots,k_s)-$equidistribution for base $\varphi$. Our definition needs two components; the first is the number of total points we require a point set to have, and the second is the number of points we require in each elementary interval. We risk repeating ourselves to emphasise the following extremely important deviation from the integer base equidistribution of digital point sets. A point set $P_m$ contained in $[0,1)^s$ of size $b^m$ is called \emph{$(k_1,\dots,k_s)-$equidistributed in base $b\in\N$} if every elementary $(k_1,\dots,k_s)-$interval in base $b$ contains exactly $b^{m-k_1-\dots-k_s}$ points from $P_m$. Intuitively, in this case the empirical distribution matches exactly with the uniform distribution on a given elementary interval. It is not immediately obvious how this idea translates when moving to an irrational base, primarily due to the fact that the number of points will clearly always be a rational number which will be impossible to match with an irrational volume of elementary interval. This means that in base $\varphi$, no matter how many points we require our equidistributed point set to have, the local discrepancy of an elementary interval (the difference between the empirical and uniform distributions) with the required number of points can never be zero. 

In line with our guiding examples, namely the van der Corput sequence in base $\varphi$ and the base $\varphi$ Hammersley point set, we will require our point set to be of size $F^m$. By the above exposition, in base $\varphi$ the local discrepancy of an elementary interval can never equal zero, and therefore we want the number of points in an elementary interval to at least \textit{minimize} its local discrepancy. In the following lemma, we see that this turns out to be a Fibonacci number of points.

\begin{lemma}\label{lem:bestapproximation}
Let $m\geq 0$ and $0\leq k\leq m+2$. The best rational approximation to $\varphi^{-k}$ with denominator $F^m$ is $F^{m-k}/F^m$ and
\[
\Bigg|\frac{F^{m-k}}{F^{m}}-\frac{1}{\varphi^k}\Bigg| = \frac{1}{\varphi^{m+2}}\frac{F^{k-2}}{F^{m}}.
\]
\end{lemma}

\begin{proof}
We apply identities \eqref{eq: fib as phi} and \eqref{eq: varphi as fib} from Section \ref{sec:definitionsandnotation}
to get
\begin{align*}
    \Bigg|\frac{F^{m-k}}{F^m}-\varphi^{-k}\Bigg| 
    &= \frac{1}{F^m}\frac{\Big|\big(\varphi^{m-k+2}-\psi^{m-k+2}\big)-\varphi^{-k}\big(\varphi^{m+2}-\psi^{m+2}\big)\Big|}{\sqrt{5}}\\
    &= \frac{1}{F^m}\frac{\big|\varphi^{-k}\psi^{m+2}-\psi^{m-k+2}\big|}{\sqrt{5}}
    =\frac{\big|\psi^{m+2}\big|}{F^m}\frac{\big|\varphi^{-k}-\psi^{-k}\big|}{\sqrt{5}}\\
    &=\frac{\big|\psi^{m+2}\big|}{F^m}\frac{\big|\varphi^{k}-\psi^{k}\big|}{\sqrt{5}}
    =\frac{1}{\varphi^{m+2}}\frac{F^{k-2}}{F^{m}}=\frac{1}{\varphi F^m+ F^{m-1}}\frac{F^{k-2}}{F^{m}}.
\end{align*}
Since this value is less than $1/F^{m}$ for all $0\leq k\leq m+2$, the best rational approximation to $\varphi^{-k}$ with denominator $F^m$ must be $F^{m-k}/F^m$.
\end{proof}

Now that we know how many points we should require in each elementary interval from the $F^m$ total points, we can finally define equidistribution in base $\varphi$.

\begin{defi}
     A point set $P_m$ contained in $[0,1)^s$ with $F^m$ points is said to be \emph{$(k_1,\dots,k_s)$-equi\-dis\-tri\-buted in base $\varphi$} if every prime elementary $(k_1,\dots,k_s)-$interval in base $\varphi$, $I$, contains exactly $F^{m-|I|}$ points from $P_m$. 
\end{defi}

It is important to point out that for a given vector $(k_1,\ldots,k_s) \in \mathbb{N}_0^s$, the value of $|I|$ is not fixed and depends on the type of prime elementary $(k_1,\ldots,k_s)-$interval being considered, as described in  \eqref{eq:logVolumeI}. One significant way in which an equidistribution in base $\varphi$ is different than an equidistribution in a natural number base, is that a $(k_1,\dots,k_s)$-equidistribution in base $\varphi$ does not imply a $(n_1,\dots,n_s)$-equidistribution in base $\varphi$ for all $0\leq n_j\leq k_j$ since although all such elementary $(n_1,\dots,n_s)$-intervals are unions of the elementary intervals in the $(k_1,\dots,k_s)$-partition, we do not require \textit{all} elementary intervals in that partition to have the correct number of points in the definition of a $(k_1,\dots,k_s)$-equidistribution. To achieve a more analogous notion of equidistribution to that of the integer base case, the following distribution property of a point set is introduced which will be useful later on.

\begin{defi}\label{def:strongequid}
A point set $P_m$ contained in $[0,1)^s$ with $F^m$ points is said to be \emph{strongly $(m_1,\dots,m_s)$-equidistributed in base $\varphi$} if every elementary interval $I$ in the  $(m_1,\dots,m_s)$-partition in base  $\varphi$, contains exactly $F^{m-|I|}$ points from $P_m$.\hjaspar{I have not defined the $(m_1,\dots,m_s)-$partition and it should be $m_j$ rather than $k_j$}
\end{defi}

The following lemma checks that indeed this strong version of equidistribution recovers the desirable property held by integer base equidistribution.

\begin{lemma}
Let $P_m$ be a point set contained in $[0,1)^s$. If $P_m$ is strongly $(m_1,\dots,m_s)$-equi\-distributed, then it is strongly $(n_1,\dots,n_s)$-equi\-distributed for all $0\leq n_j\leq m_j$. 
\end{lemma}

\begin{proof}
    It follows from Lemma \ref{lem:1d elem int characterization} that an elementary interval $I$ in the $(n_1,\dots,n_s)-$partition is the disjoint union $I=\bigcup I_i$ of elementary intervals in the $(m_1,\dots,m_s)-$partition for all $0 \leq n_j \leq m_j$. The statement will be true if the sum of the number of points in each $I_i$ under the assumption that $P_m$ is strongly $(m_1,\dots,m_s)-$equidistributed is the same as the number of points that must be in $I$ in order for $P_m$ to be strongly $(n_1,\dots,n_s)-$equidistributed. An easy consequence to  Lemmas \ref{lem:length of interval} and \ref{lem:1d elem int characterization} is that a refinement of an elementary interval $J\in[0,1)^s$ into two elementary intervals $J_1$ and $J_2$ must, without loss of generality, satisfy $|J_1| = |J|+1$ and $|J_2| = |J|+2$. The result follows by observing that 
    \[
    F^{m-|J|}=F^{m-|J|-1}+F^{m-|J|-2}=F^{m-|J_1|}+F^{m-|J_2|}
    \]
    and that each $I_i$ can be obtained from $I$ by successive refinements.
\end{proof}

\hchristiane{I adapted what is now Remark 5.8 but didn't check the rest of the paper for consistency. I think we need to make sure we want to go with the definition in Def 5.8 of tms net.}

Recall in the case of integer $(k_1, \ldots, k_s)-$equidistribution of $b^m$ total points, the number of points that should be in each elementary intervals in base $b$ is exactly $b^{m - k_1 - \cdots - k_s}$. It is then somewhat obvious upon noting that there are $b^{k_1 + \cdots + k_s}$ elementary intervals in base $b$, that we recover all $b^m$ points contained in the original set. This property is not surely guaranteed with a first glance at Definition \ref{def:strongequid}. However, the next lemma acts as a safeguard to ensure that in the definition of strong equidistribution, we indeed arrive back at $F^m$ total points when we account for the number of points in each elementary $(m_1, \ldots,m_s)-$interval in base $\varphi$. 

\begin{lemma}
    For $\mathbf{m}=(m_1, \ldots, m_s) \in \N_0^s$, let $\mathcal{P}_{\mathbf{m}}$ be an $(m_1,\dots,m_s)-$partition such that every elementary interval in it satisfies $|I|\leq m+2$. Then
    \[
    F^m=\sum_{I\in\mathcal{P}_{\mathbf{m}}}F^{m-|I|}.
    \]
    \hjaspar{Fill this in... (may remove some of those if useless). From Nathan: what was to have been filled in here?} \hnathan{Christiane and I just talked about leaving the 'parts' out and, and just having the main statement to as a check to make sure that we obtain the correct number of points in the entire net with the definition of strong equidistribution. Thoughts?}
\end{lemma}

\begin{proof}
Fix $m\geq 0$. It is clear that the $(0,\dots,0)-$partition satisfies the claim. Suppose that the $(m_1,\dots,m_s)-$partition, which we recall is denoted $\mathcal{P}_{\mathbf{m}}$, satisfies the claim and that all $I\in\mathcal{P}_{\mathbf{m}}$ satisfy $|I|\leq m+2$. Let $\mathcal{P}'_{\mathbf{m}}$ be the $(m_1,\dots,m_{j-1},m_j+1,m_{j+1},\dots,m_s)-$partition. We denote by $\mathcal I\subset \mathcal{P}_{\mathbf{m}}$ the set of elementary interval in $\mathcal{P}_{\mathbf{m}}$ that split in $\mathcal{P}'_{\mathbf{m}}$ and $\mathcal J\subset \mathcal{P}'_{\mathbf{m}}$ the set of elementary intervals in $\mathcal{P}'_{\mathbf{m}}$ that are not in $\mathcal{P}_{\mathbf{m}}$. Every $J\in\mathcal J$ is contained in exactly one $I\in\mathcal I$ and every $I\in\mathcal I$ splits as $I=J_1\cup J_2$ where $J_1,J_2\in\mathcal J$ with $|J_1|=|I|+1$ and $|J_2|=|I|+2$. Observing that $\mathcal{P}_{\mathbf{m}}\setminus \mathcal I=\mathcal{P}'_{\mathbf{m}}\setminus\mathcal J$, we calculate
\begin{align*}
F^{m}&=\sum_{I\in\mathcal{P}_{\mathbf{m}}}F^{m-|I|}
=\sum_{I\in\mathcal  I}F^{m-|I|} +\sum_{I\in \mathcal{P}_{\mathbf{m}}\setminus\mathcal I}F^{m-|I|}\\
&=\sum_{I\in \mathcal I}\big(F^{m-|I|-1}+F^{m-|I|-2}\big) +\sum_{I\in \mathcal{P}_{\mathbf{m}}\setminus\mathcal I}F^{m-|I|}\\
&=\sum_{J\in\mathcal  J}F^{m-|J|} +\sum_{J\in \mathcal{P}'_{\mathbf{m}}\setminus\mathcal J}F^{m-|J|}
=\sum_{J\in\mathcal{P}'_{\mathbf{m}}}F^{m-|J|}.
\end{align*}
\end{proof}

\begin{notation}
    For $\mathbf{m} = (m_1\dots,m_s)\in\N^s_0$, let $\supp(\mathbf m)$ be the number of non-zero coordinates of $\mathbf m$ and let 
    \begin{equation}
    \label{eq:logveck}
        \rho(\mathbf{m}) \coloneqq \sum_{j=1}^s m_j+ \supp(\mathbf{m}).
    \end{equation} 
    The purpose of this notation is to guarantee that every elementary interval $I$ in the $(m_1\dots,m_s)$-partition satisfies $|I|\leq \rho(\mathbf{m})$.
\end{notation}

An elementary interval $I$ in an $\mathbf m = (m_1,\dots,m_s)-$partition satisfies $$m_1+\dots+m_s\leq |I|\leq m_1+\dots+m_s + \supp(\mathbf m).$$ This means that when $\mathbf m$ has more non-zero coordinates, there is a greater variability in the sizes of the elementary intervals in the partition. This poses a difficulty in generalizing the definition of a $(t,m,s)-$net in an integer base because it is not clear whether the equidistribution condition should be tied to the largest elementary interval in a given partition in base $\varphi$, or the smallest. The first possibility when defining a $(t,m,s)-$net in base $\varphi$, and a more direct analogue to the definition in an integer base $b$, would require a point set to be $\mathbf k=(k_1\dots,k_s)-$equidistributed for all $k_1+\dots+k_s\leq m+2-t$. 
The second possibility would require all partitions with prime elementary intervals whose volume are all greater than or equal to $\varphi^{-(m+2-t)}$ to be equidistributed. In other words, we require the point set to be $\mathbf k=(k_1\dots,k_s)-$equidistributed for all $k_1+\dots+k_s+\supp(\mathbf k)\leq m+2-t$. 
As the point sets we will construct in Section \ref{sec:constructionspart2} work well with the second definition, we have decided to use that one. However, it may be that finding examples of higher dimensional $(t,m,s)-$nets in base $\varphi$ would make it clear which of the two definitions is best. 

\hchristiane{I would like to suggest that we use something like $\rho(\BFk)$ instead of this notation since we already have $|n|$ and $|I|$.}


\begin{defi}
\label{def:tms}
     A point set $P_m$ contained in $[0,1)^s$ with $F^m$ points is called a \emph{$(t,m,s)-$net in base $\varphi$} if it is $(k_1,\dots,k_s)$-equidistributed in base $\varphi$ for all $\mathbf k \in \mathbb{N}_0^s$ such that $\rho(\mathbf{k}) \leq m+2-t$.
\end{defi}

    

We make the following observation before introducing the notion of a $(t,s)-$sequence in base $\varphi$. A $(t,s)-$sequence in base $b\in\N$ satisfies the following: for all $n\in\N$ such that for each $0\leq i< b^m$, the least significant $m$ digits of the base $b$ expansion of $n+i$ are exactly the same as those of $i$, the set of points $\{x_n,x_{n+1},\dots,x_{n+b^m-1}\}$ is a $(t,m,s)$-net. Clearly the condition on $n$ holds if and only if $n=kb^m$ for $k,m \in \mathbb{N}_0$. For base $\MSF$, the condition on $n$ holds if and only if $n=k\odot \varphi^{m+1}$ (rather than $n=k\odot \varphi^m$).

\begin{defi}
    A sequence $\{x_0,x_1,\dots\}$ is called a \emph{$(t,s)-$sequence in base $\varphi$} if for every 
    $m \ge t$ \hchristiane{$m=t-1$ when $ t \ge 1$ would mean $\rho(\mathbf k) \le 1$ which is the same as asking $\rho(\mathbf k) \le 0$ since it cannot just be 1, so this is why I think we need to have $m \ge t$} 
    and $k \in \mathbb{N}_0$, the set $\{x_{k\odot \varphi^{m+1}},\dots x_{k\odot\varphi^{m+1} +F^m-1}\}$ is a $(t,m,s)-$net in base $\varphi$. 
\end{defi}

Analogous to the result stated as Lemma 4.22 in \cite{NIE1992} when considering integer base $b$, we have the following relationship between $(t,s)-$sequences and $(t,m,s+1)-$nets in base $\varphi$. 

\hjaspar{I do not have a reference for the integer base version of this theorem. I am probably wrong in the statement. From Nathan: maybe Wolfgang Schmid?}\hnathan{I found the reference and adapted the statement}
\begin{prop}
    If there exists a $(t,s)-$sequence in base $\varphi$, then for every $m > t $ there exists a $(t,m,s+1)-$net in base $\varphi$.
\end{prop}

\begin{proof}

Let $S = \{x_0, x_1, \ldots \}$ be a $(t,s)-$sequence in base $\varphi$. We will prove that for fixed $m \ge \max(0, t - 1)$, the point set $P$ consisting of $\left(\frac{\overline n}{\varphi^{m}}, x_{n,1},\dots,x_{n,s}\right) \in [0,1)^{s+1}$ with $0 \leq n < F^m$ is a $(t,m,s+1)-$net in base $\varphi$. 
Take a prime elementary $(k_1,\ldots,k_{s+1})-$interval in base $\varphi$,  
$$I = \prod_{j=1}^{s+1} \left[ \frac{\overline{a_j}}{\varphi^{k_j}}, \frac{\overline{a_j+1}}{\varphi^{k_j}} \right)$$ 
where $\rho(\mathbf k) \le m+2-t$. Hence the volume of $I$ is 
$\varphi^{-|I|} \ge \varphi^{-(m+2-t)}$ since
\begin{equation}\label{eq:sum}
|I| = \sum_{j=1}^{s+1} \left( k_j + |a_j| \right) \le \rho(\mathbf k) \le  m+2-t.
\end{equation}
{\bf Case 1:} $k_j>0$ for at least one $j\in \{2,\ldots,s+1\}$

We first note that in this case, $m-k_1 \ge t+1 \ge 1$.

Now $\left(\frac{\overline n}{\varphi^{m}}, x_{n,1},\dots,x_{n,s}\right) \in I$ if and only if $\overline{a_1} \varphi^{m - k_1} \leq \overline{n} < \overline{a_1+1} \varphi^{m-k_1}$, i.e. $a_1 \odot \varphi^{m-k_1} \leq n < (a_1+1) \odot \varphi^{m-k_1}$ and $$(x_{n,1},\dots,x_{n,s}) \in I' \coloneqq \prod_{j=2}^{s+1} \left[ \frac{\overline{a_j}}{\varphi^{k_j}}, \frac{\overline{a_j+1}}{\varphi^{k_j}} \right),$$ 
where $I'$ is a prime elementary interval.
Note that $\left( (a_1+1) \odot \varphi^{m-k_1} \right) - \left( a_1 \odot \varphi^{m-k_1} \right) = F^{m-k_1-|a_1|}$, i.e., a Fibonacci number of points.
Now given that the $x_n$ form a $(t,s)-$sequence in base $\varphi$ and $m-k_1-|a_1| \ge 0$ then the point set $P'$ defined as
\[
P' =\{x_n\in S:a_1 \odot \varphi^{m-k_1} \leq n < (a_1+1) \odot \varphi^{m-k_1}\}
\]
is a $(t,m-k_1-|a_1|,s)-$net in base $\varphi$. Since $I'$ is a prime elementary interval with volume $\varphi^{-|I'|}$ it contains  $F^{m-k_1-|a_1|-|I'|}=F^{m-|I|}$ points from $P'$, so $I$ contains $ F^{m-|I|}$ points from $P$ as required.

{\bf Case 2:} $k_2 = \ldots = k_{s+1}=0$ (and $k_1>0$)

In this case we only need to verify that $I = \left[ \frac{\overline{a_1}}{\varphi^{k_1}}, \frac{\overline{a_1+1}}{\varphi^{k_1}} \right)$ with $k_1 \le m+1-t$ contains $F^{m-k_1-|a_1|}$ points from $Q_m := \{\bar{n}/F^m : 0\le n<F^m\}$. 

If $m-k_1 \ge 0$, then the analysis proceeds similarly as in Case 1. That is, we note that 
$\frac{\overline n}{\varphi^{m}} \in I$ if and only if $\overline{a_1} \varphi^{m - k_1} \leq \overline{n} < \overline{a_1+1} \varphi^{m-k_1}$, i.e., $a_1 \odot \varphi^{m-k_1} \leq n < (a_1+1) \odot \varphi^{m-k_1}$. As noted in our analysis of Case 1, there are $\left( (a_1+1) \odot \varphi^{m-k_1} \right) - \left( a_1 \odot \varphi^{m-k_1} \right) = F^{m-k_1-|a_1|}$ such values of $n$, as required.

If $m-k_1 = t-1$ and $t=0$, then if $|a_1|=0$ we need to verify that $I$ contains $F^{m-|I|}=F^{-1}=1$ point from $Q_m$, while if $|a_1|=1$ then we need to verify that $I$ contains $F^{m-|I|}=F^{-2}=0$ point from $Q_m$. Since $|a_1|=0$ if and only if $I$ has a 0 coefficient for $\varphi^{-(m+1)}$ and points of the form $\bar{n}/\varphi^m$ also have a 0 coefficient for $\varphi^{-(m+1)}$, then these points can only be in intervals $I$ such that $|a_1|=0$, as required.
\end{proof}

As future research, one could derive existence results for $(t,m,s)-$nets and sequences in base $\varphi$ analogous to those that are contained in Chapter 4 of \cite{NIE1992} for nets and sequences formed from integer bases. Despite the existence of a $(0,1)-$sequence in the form of the van der Corput sequence in base $\varphi$ which will be shown shortly as Proposition \ref{prop:VDC01sequence}, the authors were unable to find a construction for a $(t,2)$-sequence for any parameter $t$, however admit that there doesn't seem to be a valid reason why one cannot exist. We were successful in deriving a two-dimensional sequence in base $\varphi$ that satisfies the following weaker definition. This construction is given in the next section and is also analysed numerically in Section \ref{sec:numericalresults}.

\begin{defi}\label{def:weaksequence}
     A sequence $\{x_0,x_1,\dots\}$ is called a \emph{weak $(t,s)-$sequence in base $\varphi$} if for all $m\geq t$, $\{x_0,x_1,\dots,x_{F^m-1}\}$ is a $(t,m,s)$-net in base $\varphi$.
\end{defi}

\section{Golden ratio point sets and sequences and their equidistribution}\label{sec:constructionspart2}

\hchristiane{Here we start using $a_j$ for the digits of $a$ but before they were used to denote numerator defining $j$th coordinate of elementary interval, e.g., in Defn 1.10. Could be confusing.}\hnathan{I think the notation here is fine. I think it is clear that these $a_j$ are digits in a digit expansion. Should a referee say something, we can come back to it.}

Now that we have redefined the concept of equidistribution for the golden ratio base, we revisit the earlier constructions from Section \ref{sec:guidingconstructions} and formally study their equidistribution properties. We start by verifying that the van der Corput sequence in base $\varphi$ is a $(0,1)$-sequence in base $\varphi$. 

\begin{remark}
Referring back to Lemma \ref{lem:bestapproximation}, we showed that given a point set containing $F^m$ points, a subinterval of $[0,1)^s$ of size $\varphi^{-k}$ should be contain $F^{m-k}$ points to minimise the local discrepancy function in that interval. For completeness, we note that the $n^{th}$ Lucas number, $L_n$ is defined by the recurrence $L_1=1$, $L_2=3$, and $L_n=L_{n-1}+L_{n-2}$ for all $n\geq 3$ and is actually the closest whole number to $\varphi^n$. Using $L_m$ rather than $F^m$ in the following definition would make sense, but the van der Corput sequence (our guiding example) would not be a $(0,1)-$sequence with respect to these numbers. 
\end{remark}

\begin{prop}\label{prop:VDC01sequence}
The van der Corput sequence in base $\varphi$ is a $(0,1)-$sequence in base $\varphi$.
\end{prop}
\begin{proof}
First we check to see that for $m \in \mathbb{N}_0$, the first $F^m$ points $\{g_0,\dots,g_{F^m-1}\}$ of the van der Corput sequence in base $\varphi$ form a $(0,m,1)-$net in base $\varphi$. Pick an arbitrary prime elementary $k$-interval $I_i(a;k)$ (recall that the adjective prime implies that $i=0,1$ only) with 
$k\leq m$ 
and $a+i=(b_{k-1}\dots b_1 b_0)_\MSF$. The volume of the selected $I_i(a;k)$ is therefore $\varphi^{-(k+i)}$, and a point $g_n$ is contained here if and only if $g_n=(.b_{k-1}\dots b_0 d_{k} \dots d_{m-1})_\varphi$. 
\hchristiane{shd this be stated as an iff statement?}\hjaspar{I think that it is iff.} 
The $n\in\N_0$ such that $g_n\in I_i(a;k)$ are all of the form
$n=(d_{m-1}\dots d_{k} b_0\dots b_{k-1})_\MSF$ or $n=j\odot\varphi^{k+i}+b$ for $0\leq j<F^{m-k-i}$, where $b=(b_0 b_1\dots b_{k-1})_\MSF$. So we conclude that there are $F^{m-k-i}$ points in the prime elementary interval $I_i(a;k)$ for $i = 0,1$ as required.

Now we also need to check the case where $k=m+1$ since $\rho(k) \le m+2$ holds iff $k \le m+1$ for nonzero $k$.
Since we know that the $(m+1)^{th}$ digit after the decimal of the first $F^m$ points is always $0$, we can say that the (prime) elementary intervals of types $I_0(a;m+1)$ and $I_1(a;m+1)$ which, from Lemma \ref{lem:length of interval}, are known to be of length $\varphi^{-m-1}$ and $\varphi^{-m-2}$, have one point and zero points in them respectively. This is exactly what is required to be $(m+1)-$equidistributed. Thus the first $F^m$ points of the van der Corput sequence are a $(0,m,1)-$net in base $\varphi$. It might be of interest for the reader to refer back to Lemma \ref{lem:1d elem int characterization} to clarify how intervals in $\mathcal{P}_m$ are refined in $\mathcal{P}_{m+1}$. 

The $(0,1)-$sequence property follows directly from the fact that the first $m+1$ digits of the points $\{g_{k\odot \varphi^{m+1}},\dots g_{k\odot\varphi^{m+1} +F^m-1}\}$ are exactly the same as the first $m+1$ digits of the points $\{g_0,\dots,g_{F^{m-1}-1}\}$. Since the latter set is a $(0,m,1)$-net in base $\varphi$, the former must be too.
\end{proof}

We now show an equivalent two-dimensional result regarding the distribution of the Hammersley point set in base $\varphi$. The notation \eqref{not:arbprimeinterval} will be employed in this proof and therefore urge the reader to revise the same.

\begin{prop}\label{prop:Ham0m2net}
The Hammersley set $H_m$ in base $\varphi$ with $F^m$ points is a $(0,m,2)$-net in base $\varphi$.
\end{prop}

\begin{proof}
\hchristiane{I think we need to be more clear as to how this is derived from Defn 1.10 and Defn 1.13}
We recall from Definition \ref{def:tms} that a point set is called a $(0,m,2)-$net in base $\varphi$ if it is $(k_1, k_2)-$equi\-dis\-tri\-buted in base $\varphi$ for all $\rho(\mathbf{k}) \leq m+2$. Suppose that $I(a,b;k_1,k_2)$ is a two-dimensional prime elementary $(k_1,k_2)-$interval with $\rho(\mathbf{k}) \leq m + 2$, $a=(a_{k_1-1}\dots a_1a_0)_\MSF$ and $b=(b_{k_2-1}\dots b_1b_0)_\MSF$. As a key observation, note that the points from $H_m$ that belong in $I(a,b;k_1,k_2)$ are of the form
\[
\big((.a_{k_1-1}\dots a_0 d_{k_1}\dots d_{m-k_2-1} b_0\dots b_{k_2-1})_\varphi,
(.b_{k_2-1}\dots b_0 d_{m-k_2-1}\dots d_{k_1} a_0\dots a_{k_1-1})_\varphi\big).
\]

We have two cases to consider. Let's firstly suppose that both $k_1, k_2 > 0$, i.e., our elementary $(k_1,k_2)-$interval is such that $k_1 + k_2 \leq m$. If $|a|=0$ then $d_{k_1}$ may be either $0$ or $1$, but if $|a|=1$, then $d_{k_1}$ must be $0$. Similarly if $|b|=0$ then $d_{m-k_2-1}$ may be either $0$ or $1$, but if $|b|=1$ then $d_{m-k_2-1}$ must be $0$. Since $(d_{k_1}\dots d_{m-k_2-1})_\MSF$ must be a binary string with no consecutive ones and incorporating the possibility that $d_{k_1}$ or $d_{m-k_2-1}$ must be zero, we can conclude that there are exactly $F^{m-k_1-k_2-|a|-|b|} = F^{m-|I|}$\hjaspar{maybe we should have a lemma counting the number of binary strings with no consecutive ones in the previous section} \hchristiane{yes} points in the elementary interval $I(a,b;k_1,k_2)$. \hchristiane{I also note that this corresponds to $F^{m-|I|}$ with $|I|$ from Defn 1.10, because the ``types'' $i$ and $j$ are equal to the last digit of $k_1$ and $k_2$}  
\hchristiane{This is describing the case where exactly one of $k_1$ or $k_2$ is 0. Would help if we started at the beginning and explain we have these two cases to consider.}

As the second and final case, we now assume that one of $k_1$ or $k_2$ is equal to $0$, i.e., our prime elementary interval is such that $\rho(\mathbf{k}) = k_1+k_2 + 1\leq m+2$ or, $k_1+k_2 \leq m+1$. We can deduce from the previous proposition that $H_m$ is both $(k_1,0)-$equidistributed and $(0,k_2)-$equidistributed because the first $F^m$ points of the van der Corput sequence are $k-$equidistributed for $0 \le k \le m+1$. Thus, $H_m$ is a $(0,m,2)-$net in base $\varphi$.
\end{proof}

\subsection{Two-dimensional sequences in base $\varphi$}

\hnathan{Place the discussion of creating a $(t,2)-$sequence here.}
\hjaspar{I can say what I tried to create such sequences, but I don't have anything rigorous.}

As mentioned at the end of Section \ref{sec:equid}, the authors were unable to find a $(t,2)-$sequence in base $\varphi$ for any quality parameter $t$. We were, however, successful in obtaining a so-called weak $(1,2)$-sequence in base $\varphi$ as per Definition \ref{def:weaksequence}. We present two important lemmas which are subsequently used to describe how a weak two-dimensional sequence in base $\varphi$ is obtained. To be used as an aid when reading Lemmas \ref{lem: point distribution} and \ref{lem: Extending points}, Figure \ref{fig:22 32 23 partitions} may be insightful to understand how two-dimensional partitions in base $\varphi$ are refined.

\begin{lemma}\label{lem: point distribution}
Let $P_m$ be a $(1,m,2)$-net in base $\varphi$ and let $1\leq k < m$. Then each group of four prime elementary intervals
\[
I_{i,j}(a,b;m-k,k)=\Bigg[\frac{\overline {a+i}}{\varphi^{m-k}},\frac{\overline{a+i+1}}{\varphi^{m-k}}\Bigg)\times\Bigg[\frac{\overline {b+j}}{\varphi^{k}},\frac{\overline{b+j+1}}{\varphi^{k}}\Bigg),\quad i,j\in\{0,1\}
\]
together contain exactly three points from $P_m$ distributed amongst them in one of two ways:
\begin{enumerate}
    \item either there are two points in $I_{0,0}(a,b;m-k,k)$ and one point in $I_{1,1}(a,b;m-k,k)$, or
    \item there is one point in each of $I_{0,0}(a,b;m-k,k)$, $I_{1,0}(a,b;m-k,k)$, and $I_{0,1}(a,b;m-k,k)$.

\end{enumerate}
Moreover, elementary intervals $I_{i,j}(a,b;m-k,k)$ in the $(m-k,k)$-partition with at least one of $i$ or $j$ equal to $2$ have exactly the right number of points from $P_m$ (either $2$ or $1$ depending on their size).\hjaspar{maybe I should be specific.}
\end{lemma}

\begin{figure}[t!]
      \centering
      \subfigure[]{\includegraphics[width=0.35\textwidth]{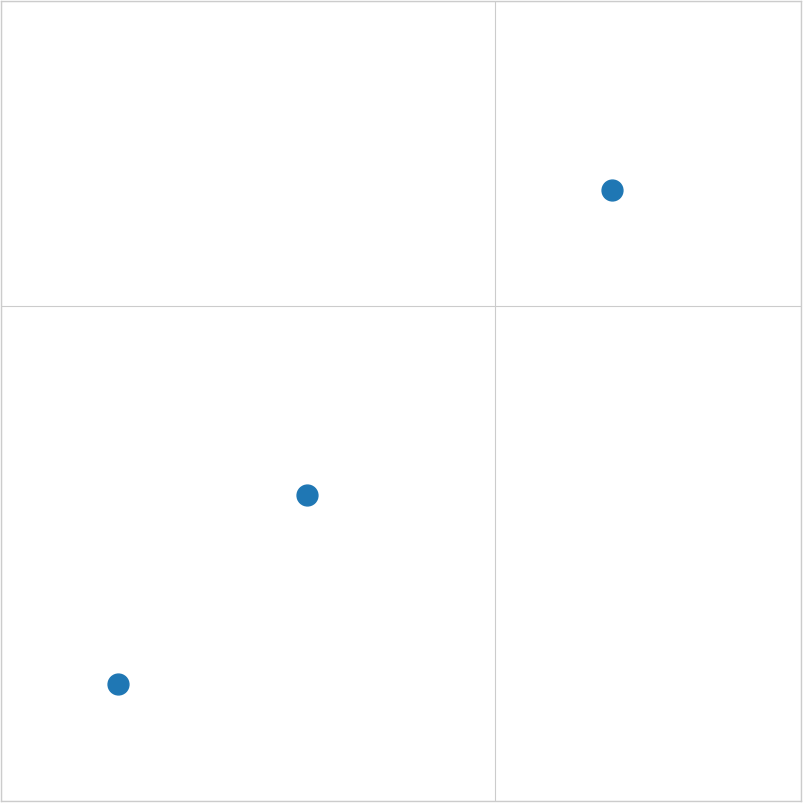}}
      \hspace{5mm}
      \subfigure[]{\includegraphics[width=0.35\textwidth]{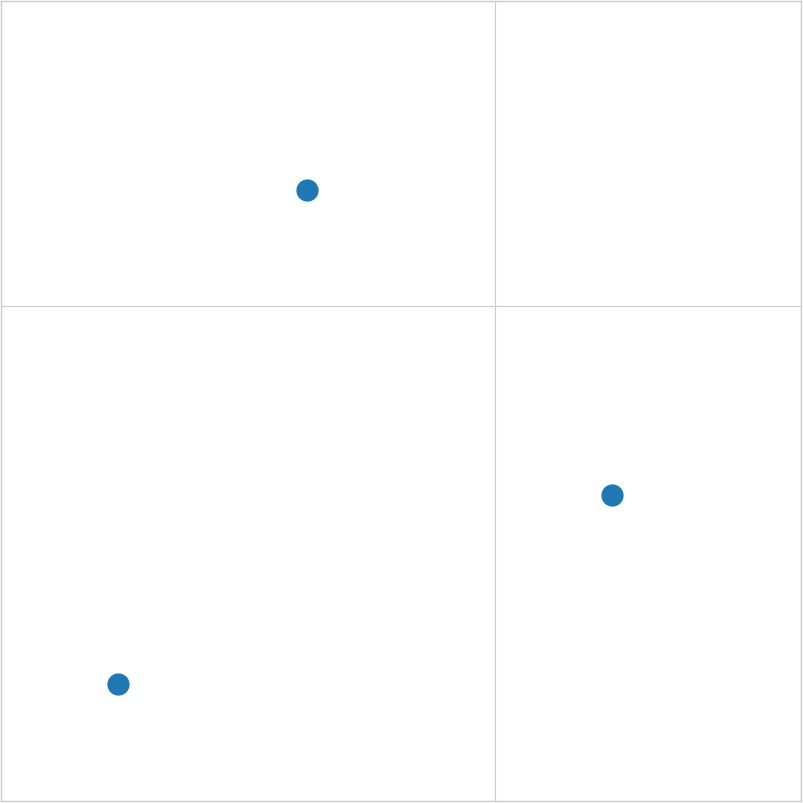}}
      \caption{The two possible point arrangements from Lemma \ref{lem: point distribution} in each group of four elementary $(m-k,k)$-intervals $I_{i,j}(a,b;m-k,k)$ where $a$ and $b$ are fixed and $i,j\in\{0,1\}$.}
      \label{fig:point distribution}
 \end{figure}

\begin{proof}
Suppose $P_m$ is a $(1,m,2)$-net in base $\varphi$. We see from Lemma \ref{lem:1d elem int characterization} that 
\begin{align*}
I_{0,0}(a\odot\varphi^{-1},b;m-k-1,k) &= I_{0,0}(a,b;m-k,k)\cup I_{1,0}(a,b;m-k,k),\\
I_{0,0}(a,b\odot\varphi^{-1};m-k,k-1) &= I_{0,0}(a,b;m-k,k)\cup I_{0,1}(a,b;m-k,k),\\
I_{0,1}(a\odot\varphi^{-1},b;m-k-1,k) &= I_{0,1}(a,b;m-k,k)\cup I_{1,1}(a,b;m-k,k), \text{ and}\\
I_{1,0}(a,b\odot\varphi^{-1};m-k,k-1) &= I_{1,0}(a,b;m-k,k)\cup I_{1,1}(a,b;m-k,k)
\end{align*}
and that these are prime elementary $(m-k-1,k)$-intervals or prime elementary $(m-k,k-1)$-intervals. By the definition of $(t,m,s)$-net in base $\varphi$, the first two intervals each must have exactly two points from $P_m$ and the second two must each contain exactly one point from $P_m$. If $I_{1,1}(a,b;m-k,k)$ does not contain a point from $P_m$, then, as each of $I_{0,1}(a\odot\varphi^{-1},b;m-k-1,k)$ and $I_{1,0}(a,b\odot\varphi^{-1};m-k,k-1)$ contain exactly one point from $P_m$, both $I_{0,1}(a,b;m-k,k)$ and $I_{1,0}(a,b;m-k,k)$ must contain exactly one point from $P_m$ and since $I_{0,0}(a\odot\varphi^{-1},b;m-k-1,k)$ contains exactly two points from $P_m$, $I_{0,0}(a,b;m-k,k)$ must contain exactly one point from $P_m$ (this is the right point set in Figure \ref{fig:point distribution}). Similar reasoning shows that when $I_{1,1}(a,b;m-k,k)$ contains exactly one point from $P_m$, there must be exactly two points from $P_m$ in $I_{0,0}(a,b;m-k,k)$ (this is the left point set in Figure \ref{fig:point distribution}). As $I_{1,1}(a,b;m-k,k)$ cannot contain more than one point from $P_m$ since then $I_{0,1}(a\odot\varphi^{-1},b;m-k-1,k)$ would have too many points, these are the only two cases. Again by Lemma \ref{lem:1d elem int characterization}, for $i,j\in\{0,1\}$,
\begin{align*}
I_{i,2}(a,b;m-k,k)&=I_{i,1}(a,b\odot\varphi^{-1};m-k,k-1),\\
I_{2,j}(a,b;m-k,k) &= I_{1,j}(a\odot\varphi^{-1},b;m-k-1,k), \text{ and}\\
I_{2,2}(a,b;m-k,k) &= I_{1,1}(a\odot\varphi^{-1},b\odot\varphi^{-1};m-k-1,k-1).
\end{align*}
From the right-hand side we see all three are prime
elementary $(k_1,k_2)$-intervals with $\rho(k_1,k_2) \le m+1$,  which means that they must contain the correct number of points from $P_m$, which is what the last part of the statement of the lemma claims. 
\end{proof}

\begin{lemma}
    Let $P_m$ be a point set in $[0,1)^2$ with $F^m$ points. If $P_m$ is strongly $(m-k,k)-$equi\-distributed for $0 \le k \le m$, then $P_m$ is a $(1,m,2)-$net in base $\varphi$.
\end{lemma}

\begin{proof}
The strong $(m-k,k)$-equidistribution implies all elementary intervals in a $(k_1,k_2)$-partition with $k_1+k_2 \le m$ and $k_1,k_2 \ge 0$ have the correct number of points. Now consider $k_1,k_2 \ge 1$ such that  $\rho(k_1,k_2) \le m+1$. Therefore $k_1 + k_2 \le m-1$ and thus all (prime) elementary intervals in the $(k_1,k_2)-$partition have the correct number of points. If only one of $k_1$ or $k_2$ is non-zero and equal to $k$ with $k\le m$, then $\rho(k_1,k_2) = k+1 \le m+1$. By our assumption on $P_m$ all (prime) elementary intervals in the corresponding $(0,k)$ or $(k,0)$ partition have the correct number of points.
\end{proof}

\hchristiane{If we ask for strong $(m+1-k,k)$-equidistribution, then we get a $(0,m,2)-$net. This means checking that in $(m+1-k,k)-$ partitions, we have one point in the interval of type $I_{0,0}$ and no point in all the others. In turn, this seems to only request that points are in the lower left corner in the elementary intervals of the $(m-k,k)$-partition.}

\begin{lemma}\label{lem: Extending points}
Let $S_m$ be a $(1,m,2)$-net in base $\varphi$ that is strongly $(m+1-k,k)$-equidistributed for $0 \le k \le m+1$. There exists a $(0,m-1,2)$-net, $P_{m-1}$ in base $\varphi$ that makes $S_m\cup P_{m-1}$ a $(1,m+1,2)$-net in base $\varphi$ that is strongly $(m+2,0)$-equidistributed and strongly $(0,m+2)$-equidistributed.
\end{lemma}

\begin{proof}
The first step is to show that the first coordinate of the points in $P_{m-1}$ are uniquely determined by the requirement that $S_{m}\cup P_{m-1}$ is strongly $(m+1,0)$-equidistributed \hjaspar{change this}. The $(m,0)$-partition contains three types of elementary intervals, namely $I_{0,0}(a,0;m,0)$, $I_{1,0}(a,0;m,0)$ and $I_{2,0}(a,0;m,0)$, each of which contains exactly one point from $S_m$. Both $I_{0,0}(a,0;m,0)$ and $I_{2,0}(a,0;m,0)$ split in the $(m+1,0)$-partition as $I_{0,0}(\tilde{a},0;m+1,0)\cup I_{1,0}(\tilde{a},0;m+1,0)$ (where $\tilde{a}$ is either $a\odot\varphi$ or $a\odot\varphi+3$), with each pair having one point from $S_m$ between them, and $I_{1,0}(a,0;m,0)=I_{2,0}(a\odot\varphi,0;m+1,0)$ has one point from $S_m$ in it. Each of the three types of elementary intervals in the $(m+1,0)$-partition needs exactly one point from $S_m\cup P_{m-1}$ for $S_m\cup P_{m-1}$ to be strongly $(m+1,0)$-equidistributed (since type 0 and type 2 need $F^{m+1-(m+1)}$ and type 1 needs $F^{m+1-(m+2)} = F^{-1}=1$), thus the points in $P_{m-1}$ must be in those $I_{i,0}(a,0;m+1,0)$ that do not contain a point from $S_m$. We define ${}_0P$ to be the point set 
\[
{}_0P=\Bigg\{\Bigg(\frac{a+i}{\varphi^{m+1}},0\Bigg)\in[0,1)^2: I_{i,0}^{m+1,0}(a,0)\cap S_m=\emptyset, 0\leq a=(d_{m}\dots d_2 0 0)_\MSF<F^{m+1}, i\in\{0,1\}\Bigg\}.
\]
We remark for clarity, that there are $F^{m-1}$ points in ${}_0P$ despite the fact that these points contain $m+1$ digits. 

It is then clear that $S_m\cup {}_0P$ is strongly $(m+1,0)$-equidistributed. Also, since for valid $a$ and $i\in\{0,2\} $
\[
I_{i,0}(a,0;m,0) = I_{0,0}(\tilde{a},0;m+1,0)\cup I_{1,0}(\tilde{a},0;m+1,0)
\]
contains exactly one point from ${}_0P$ and the points from ${}_0P$ only appear in those elementary intervals, ${}_0P$ is strongly $(m,0)$-equidistributed. Thus, we have found the first coordinate of the points in ${}_0P$. 

We find the digits of the second coordinate 
by induction with ${}_kP$ denoting the point set which has had digit $1$ up to $k$ of the second coordinate chosen.  The induction hypothesis is that ${}_{k-1}P$ is strongly $(m-n,n)$-equidistributed for all $0\leq n<k$, each point in ${}_{k-1}P$ has zeros after the $(k-1)^{th}$ digit in the second coordinate, and $S_m\cup {}_{k-1}P$ is strongly $(m+1-n,n)$-equidistributed for all $0\leq n<k$. We have already done the base case.\hjaspar{I have only shown that that $P_k$ satisfies the $n=k$ part of the induction hypothesis.}\hchristiane{Not sure what that means}

\hchristiane{I stopped here. I need to go through this part of the proof to validate the assumptions in the proposition and the I.H. To do so I need to better understand the connections between the net property and the $(m-k,k)-$s.e.\ type property.}
Suppose that ${}_{k-1}P$ satisfies the induction hypothesis. We begin with the groups of four prime elementary $(m+1-k,k)$-intervals, $I_{i,j}(a,b;m+1-k,k)$ where $i,j\in\{0,1\}$. Due to Lemma \ref{lem: point distribution} and the fact that $S_m$ is a $(1,m,2)-$net, in Figure \ref{fig:Cases for adding points} we see the two possible cases of how the points from $S_m$ (shown in gray) can be distributed amongst these four elementary intervals, with the first case being the top row and the second case being the bottom row. By Lemma \ref{lem:1d elem int characterization}, the union of the top and bottom elementary $(m+1-k,k)$-intervals on the left-hand side is an elementary interval in the $(m+1-k,k-1)$-partition, 
\[
I_{0,j}(a,b\odot\varphi^{-1};m+1-k,k-1) = I_{0,0}(a,b;m+1-k,k)\cup I_{0,1}(a,b;m+1-k,k),
\]
where $j$ is either $0$ or $2$. Since ${}_{k-1}P$ is strongly $(m+1-k,k-1)$-distributed the elementary interval on the left hand side has one point from ${}_{k-1}P$, which must be on the bottom edge because only the initial $k-1$ digits of the second coordinate of this point can be non-zero. Similarly the union of the top and bottom elementary intervals on the right-hand side is also an elementary interval in the $(m+1-k,k-1)$-partition
\[
I_{1,j}(a,b\odot\varphi^{-1};m+1-k,k-1) = I_{1,0}(a,b;m+1-k,k)\cup I_{1,1}(a,b;m+1-k,k),
\]
where $j$ is either $0$ or $2$. Again, the elementary interval on the left-hand side must have exactly one point from ${}_{k-1}P$ located on the bottom edge. Thus, there is a point from ${}_{k-1}P$ on the bottom edge of $I_{i,0}(a,b;m+1-k,k)$ for $i\in\{0,1\}$, on the left-hand side of Figure \ref{fig:Cases for adding points} we see the two possible distributions of points from $S_{m}\cup {}_{k-1}P$ in a group of four elementary intervals with the top being Case 1 and the bottom being case 2 from Lemma \ref{lem: point distribution}. If we are in case 1, then the left point from ${}_{k-1}P$ will have its $k^{th}$ digit changed to a $1$ in ${}_kP$ while the right point will have its $k^{th}$ digit equal to $0$ and the reverse will happen in Case 2. In this way, all four prime elementary intervals have the right number of points from $S_m\cup {}_kP$. 

\begin{figure}[t!]
     \centering
     \includegraphics[scale = .25]{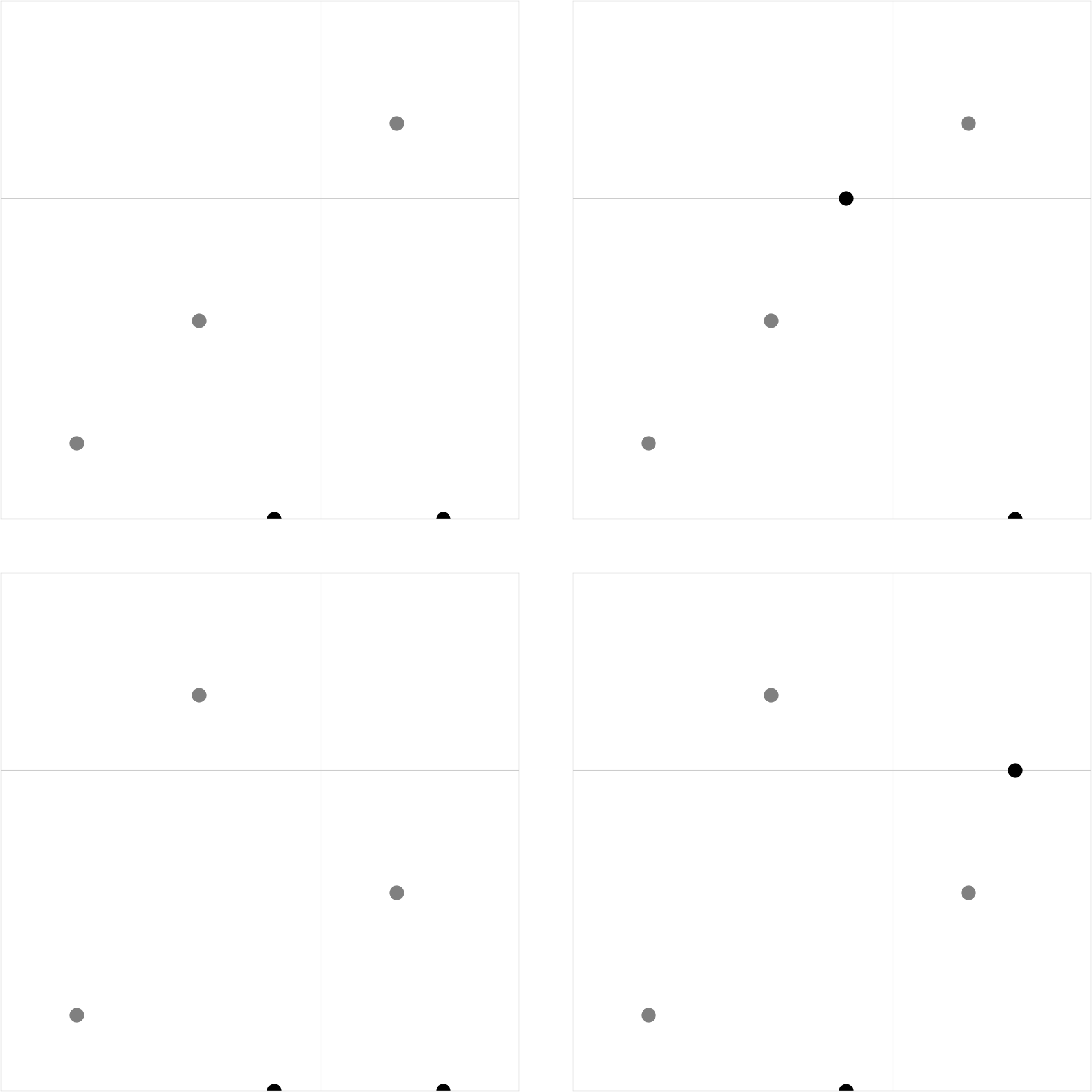}
     \caption{The cases for Lemma \ref{lem: Extending points}}\label{fig:Cases for adding points}
\end{figure}

Now, the union of the bottom two and top two elementary intervals in the $(m+1-k,k)$-partition become elementary intervals in the $(m-k,k)$-partition,
\begin{align*}
    I_{i,0}(a\odot\varphi^{-1},b;m-k,k)&=I_{0,0}(a,b;m+1-k,k)\cup I_{1,0}(a,b;m+1-k,k)\text{ and }\\
     I_{i,1}(a\odot\varphi^{-1},b;m-k,k)&=I_{0,1}(a,b;m+1-k,k)\cup I_{1,1}(a,b;m+1-k,k),
\end{align*}
where $i$ is either $0$ or $2$. It is clear that in both cases, there is one point from ${}_{k}P$ in each of the two elementary intervals, which is how many there needs be in order for ${}_kP$ to be strongly $(m-k,k)$-equidistributed. We note that every elementary interval $I_{i,j}(a,b;m-k,k)$ with $i\in\{0,2\}$ and $j\in\{0,1\}$ splits in the $(m+1-k,k)$-partition in this way which means every elementary interval in the $(m-k,k)$-partition of these types has the right number of points from ${}_{k}P$.

Next we look at the pairs of elementary intervals $I_{2,j}(a,b;m-k,k)$ where $j\in\{0,1\}$. Both of these are elementary $(m-k,k)$-intervals and thus each has exactly one point from $S_m$. As before, these two elementary intervals combine to form the elementary interval $I_{2,0}(a,b\odot\varphi^{-1};m-k,k-1)$ which must have one point from ${}_{k-1}P$ located on its lower edge, thus $I_{2,0}(a,b;m-k,k)$ contains a point from ${}_{k-1}P$ and $I_{2,1}(a,b;m-k,k)$ contains no points from ${}_{k-1}P$. Thus $I_{2,0}(a,b;m-k,k)$ contains two points from $S_m\cup {}_{k-1}P$ and $I_{1,0}(a,b;m-k,k)$ contains one point from $S_m\cup {}_{k-1}P$ which is how many they should have. This means that $k^{th}$ digit of every point from ${}_{k-1}P$ in an elementary interval of type $I_{2,0}(a,b;m-k,k)$ remains $0$ in ${}_kP$.

\begin{figure}[t!]
    \centering
    \subfigure[]{\includegraphics[width=0.48\textwidth]{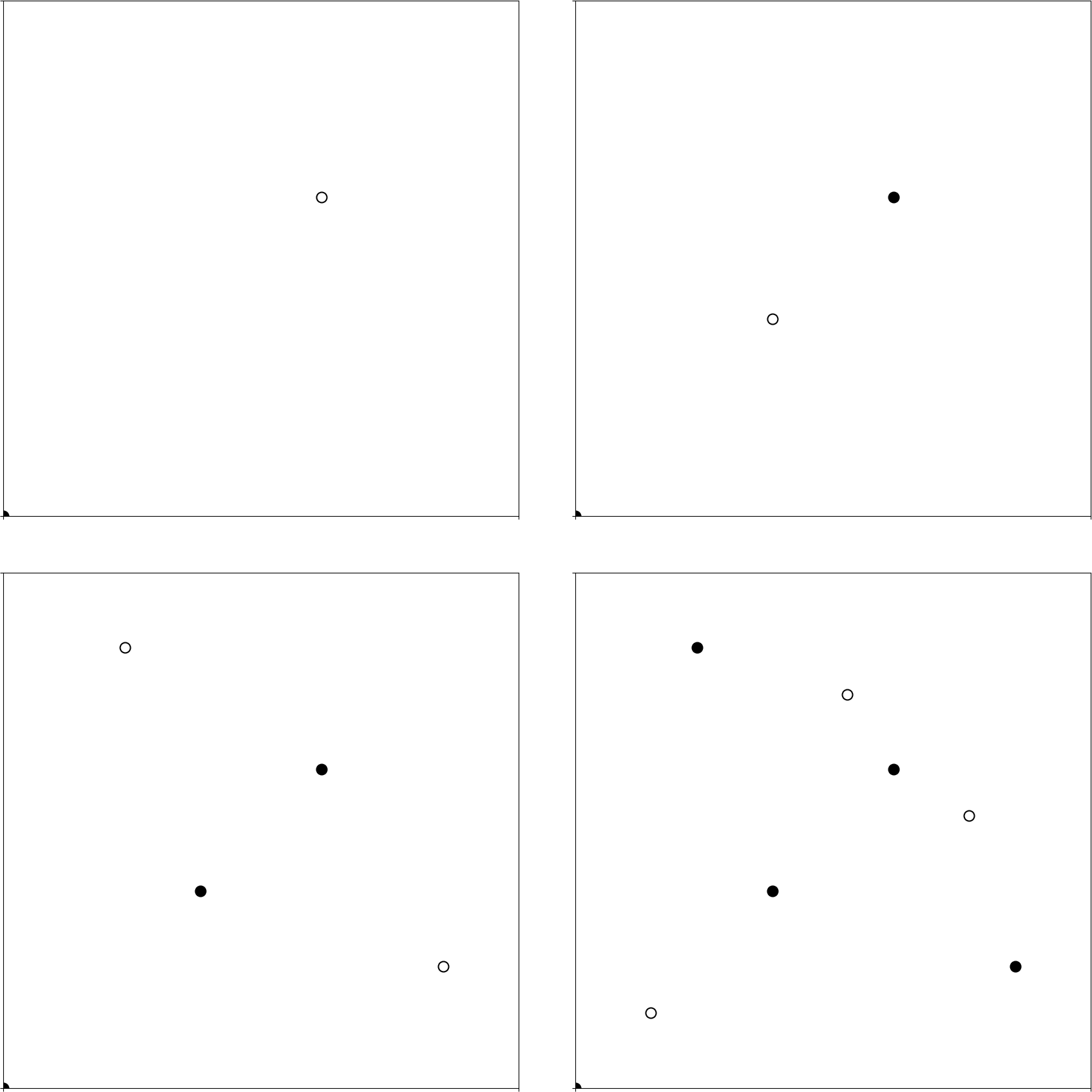}} 
    \hspace{3mm}
    \subfigure[]{\includegraphics[width=0.48\textwidth]{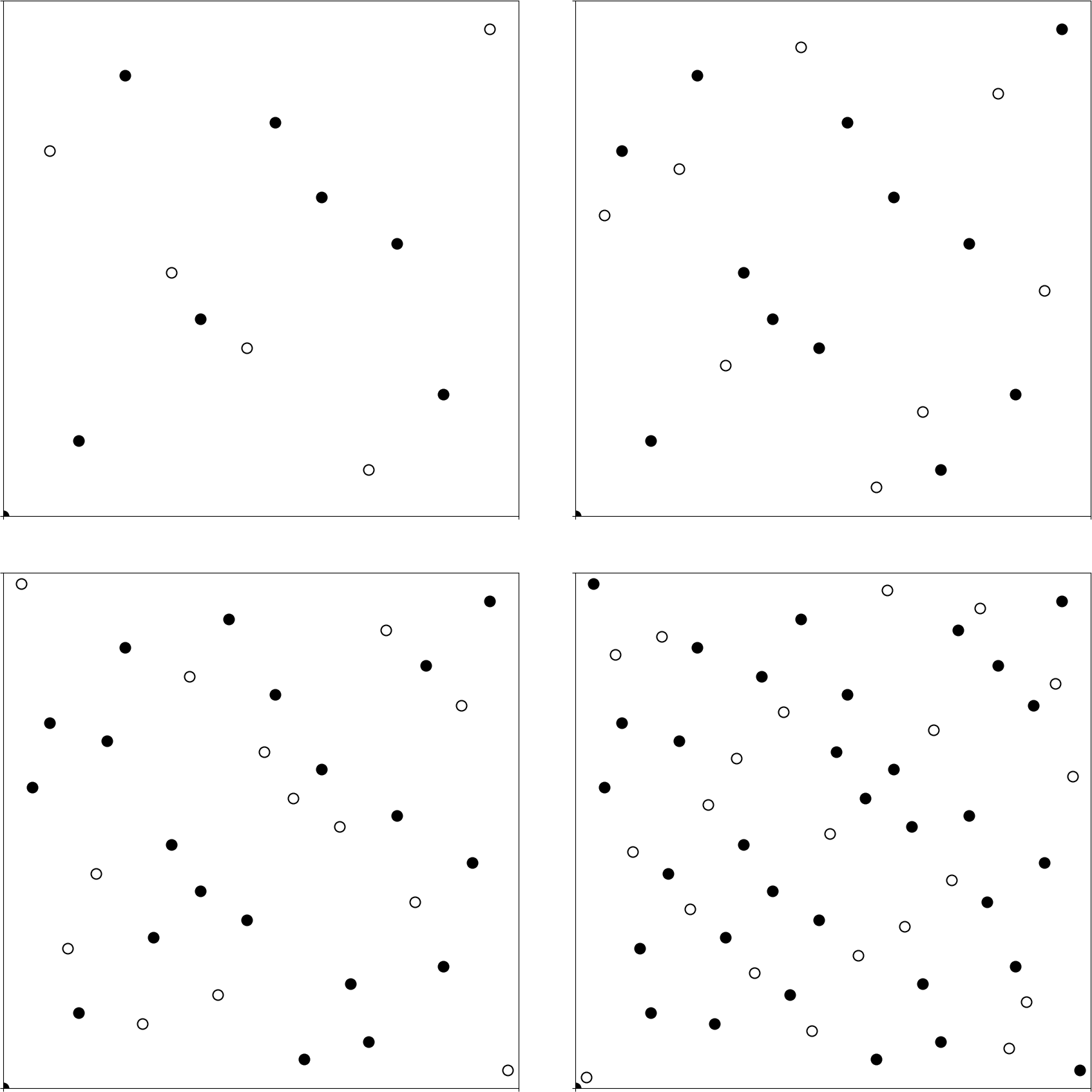}} 
    \caption{The sequence generated by applying Lemma \ref{lem: Extending points} repeatedly with the origin as the starting point. The points to be added at each step are circles. The first four steps are given in (a) and the second four in (b).}
    \label{fig:extensions from lemma steps 1 to 8}
\end{figure}

The elementary intervals in the $(m-k,k)$ partition of the form  $I_{1,0}(a,b;m-k,k)$ appear in the $(m+1-k,k)$-partition as $I_{2,0}(a\odot\varphi,b;m+1-k,k)$ which, we have seen, have one point from $P_k$ while $I_{1,1}(a,b;m-k,k)=I_{2,1}(a,b;m+1-k,k)$ has no point from ${}_kP$. Thus all elementary intervals of this type have the correct number of points from ${}_kP$ that they need to have in order for ${}_kP$ to be strongly $(m-k,k)$-equidistributed.

We turn our attention to pairs of elementary intervals of the form $I_{i,2}(a,b;m+1-k,k)$ where $i\in\{0,1\}$. Since these elementary intervals appear in the $(m+1-k,k-1)$-partition as $I_{i,2}(a,b;m+1-k,k)=I_{i,1}(a,b\odot\varphi^{-1};m+1-k,k-1)$ they both have one point from $S_m$ with the one with $i=0$ having one point from ${}_{k-1}P$ and the one with $i=1$ having no points from ${}_{k-1}P$. This means that these elementary intervals have the right number of points from $S_m\cup {}_{k-1}P$ and we need not change the $k^{th}$ digit of any point from ${}_{k-1}P$ in elementary intervals of type $I_{0,2}(a,b;m+1-k,k)$. Elementary intervals in the $(m+1-k,k)$-partition of the form $I_{0,2}(a\odot\varphi^{-1},b;m+1-k,k)$ are all the union of these pairs of intervals in the $(m+1-k,k-1)$-partition and therefore contain one point from ${}_kP$, which is how many they need in order for ${}_kP$ to be strongly $(m+1-k,k)$-equidistributed.

Finally, elementary intervals of the form 
\[
I_{2,2}(a,b;m+1-k,k)=I_{2,1}(a,b\odot \varphi^{-1};m+1-k,k-1)=I_{1,2}(a\odot\varphi^{-1},b;m-k-1,k)
\]
have one point from $S_m$ and one point from ${}_{k-1}P$ and the point from ${}_{k-1}P$ goes into ${}_{k}P$ with its $k^{th}$ digit equal to $0$, giving us what we want. We have shown that ${}_{k}P$ satisfies the induction hypothesis.

We now have an $m$-digit point set ${}_mP$ that satisfies the induction hypothesis. The last step is to find the $(m+1)^{th}$ digit of the second coordinate of each point in $P_{m-1}$. To do this we simply choose that digit in such a way as to make $S_m\cup P_{m-1}$ strongly $(0,m+1)-$equidistributed. This is nearly identical to how we chose the first coordinate of $P_{m-1}$ which was possible because $S_m$ is strongly $(m+1,0)-$equidistributed.
\end{proof}

With the work completed in this subsection, we have the following existence result.

\begin{prop}
    There exist weak $(1,2)$-sequences in base $\varphi$.
\end{prop}

\begin{proof}
    There are many two-dimensional point sets that satisfy the assumptions of Lemma \ref{lem: Extending points}. Indeed, take $S_0=\{ x_0\}$ which is a $(0,0,2)-$net in base $\varphi$ (note that we in fact only need a $(1,0,2)$-net) and is both $(1,0)-$ and $(0,1)-$equidistributed in base $\varphi$ whenever $x_0\in I_{0,0}(0,0;1,1)$. Applying Lemma \ref{lem: Extending points} to $S_0$ yields a point set $\{x_0, x_1\}$ that satisfies the assumptions of that same lemma. In particular, $\{x_0, x_1\}$ is a $(1,1,2)$-net in base $\varphi$. Continuing by induction, we get a sequence of points $\{ x_0, x_1, x_2,\dots\}$ such that each $\{x_0,\dots,x_{F^m-1}\}$ is a $(1,m,2)$-net in base $\varphi$ and thus satisfying the definition of a weak $(1,2)$-sequence where each additional $F^{m-1}$ points $\{x_{F^{m}},\dots,x_{F^{m+1}-1}\}$ are obtained by applying Lemma \ref{lem: Extending points} to $\{x_0,\dots,x_{F^m-1}\}$. 
\end{proof}

\subsection{Scrambling in base $\varphi$} The notion of \textit{scrambling} was first introduced by A. Owen \cite{OWEN1995} with the intention of introducing an element of randomization to digital point sets and sequences. We recall the basic notion of Owen scrambling in one dimension for the reader's benefit and in order to discuss the attempt of generalization to base $\varphi$.

In base $b\geq 2$, one selects a permutation $\pi$ which is uniformly distributed over the $b!$ permutations of $\{0,\ldots,b-1\}$. Then, the first digit in the base $b$ expansion of some $x\in [0,1)$ is permuted by $\pi$, and subsequently the second digit is permuted by another random permutation that depends on the  first digit. Continuing in a similar manner, the permutation applied to the $k^{th}$ digit depends on the values of the preceding $k-1$ digits. Owen's scrambling possesses two very useful features; firstly, the randomization of points preserves the equidistribution properties of $(t,m,s)-$nets and $(t,s)-$sequences in base $b$, i.e. the $t$ parameter is preserved and secondly, each of the scrambled points are uniform in $[0,1)$; for the details and proofs of these claims, see \cite{OWEN1995}. It is therefore important that scrambling in base $\varphi$ should also possess these desirable properties.

In an attempt to generalize scrambling to the golden ratio base, let us take a point set $P$ and take the $i^{th}$ point denoted by $a_i = (.a_0a_1 \ldots a_{m-1} a_m a_{m+1} \ldots )_{\varphi}.$ At any given stage of the scrambling, we want to permute the $m^{th}$ digit of the point depending on the previous $(m-1)$ digits. However note that due to the restraint in base $\varphi$ that we are not allowed consecutive ones in a digit expansion and consequently, we can only freely permute the $m^{th}$ digit $a_m$ if the two neighboring digits $a_{m-1}$ and $a_{m+1}$ have both been set to 0 after scrambling. Unfortunately, this turns out to be a significant hurdle for scrambling nets and sequences in base $\varphi$ since there exists a dependence on successive permutations that is not present in the integer base $b$ scrambling.

Thus, should a valid scrambling algorithm exist, it would need to be clever in construction to overcome the problem caused due to the restriction on digits in base $\varphi$. We consider investigating a method for scrambling irrational based nets and sequences as a future project. As a first step, utilising the notion of \textit{quasi-equidistribution in base $b$} first introduced in \cite{WIARTLEMIEUXDONG2021}, numerical investigation by the authors of this paper (not contained in this text) suggests the following result.

\begin{conjecture}\label{con:quasi-equidistribution}
    The Hammersley point set in base $\varphi$ is \textit{completely quasi-equidistributed} in base $2$.
\end{conjecture}

Should the above conjecture hold true, this could suggest that one can scramble the base $\varphi$ Hammersley set in base $2$ and still retain desirable uniformity in the newly randomized point set.

\section{Numerical calculations of discrepancy}\label{sec:numericalresults}

\begin{figure}[t!]
     \centering
     \includegraphics[width=16cm]{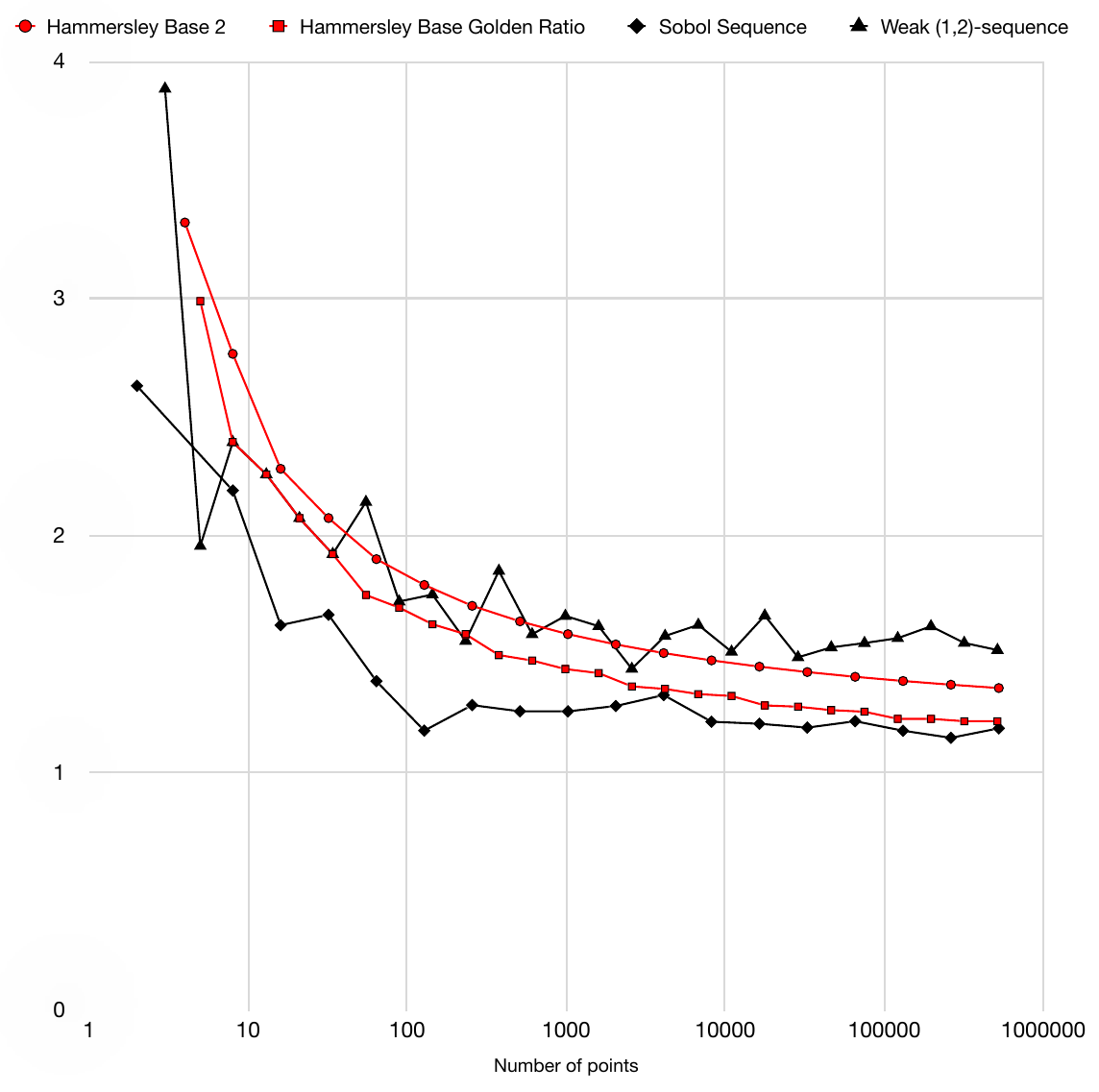}
     \caption{The star discrepancy $D_N^{*}(\cdot)$ of the Hammersley point set in base $2$ and $\varphi$, the weak $(1,2)-$sequence and the Sobol' sequence in two-dimensions normalized by $\log N/N$.}\label{fig:star disc}
\end{figure}

In this section, we numerically evaluate the performance of the base $\varphi$ point sets and sequences that have been constructed so far in this paper. We focus on the two-dimensional constructions and give extensive numerical results to evaluate their performance with respect to the discrepancy. In particular, we evaluate the star discrepancy (and $L_2-$discrepancy; see Appendix \ref{app:additional figures - discrepancy}, Figure \ref{fig:L2 star disc}) for four point sets and sequences: the Hammersley construction in base $2$ and base $\varphi$, the weak $(1,2)-$sequence obtained from Section \ref{sec:constructionspart2} and lastly, the classical Sobol' sequence in dimension two. The discrepancy values in Figure \ref{fig:star disc} have been normalized by $\log N/N$, and the evident convergence of all four point sets means that the vertical axis is actually tracking the constant factor independent of $N$ for the star discrepancy.  

Studying Figure \ref{fig:star disc}, we can make several observations. Firstly, the discrepancy of the Hammersley construction in base $\varphi$ clearly outperforms the classical base $2$ construction for the entire range of values of $N$. Therefore we can conclude that rather than restricting oneself to the natural numbers, one can improve distribution properties of the Hammersley point set construction by carefully choosing a base from the set of positive reals. This naturally leads to our first open question for further research:

\vspace{0.2cm}
\noindent
\textbf{Question 1.} What choice of base $\gamma \in \mathbb{R}$ yields the smallest star discrepancy of the Hammersley point set in base $\gamma$?

\vspace{0.2cm}
In the next and final section, we begin to explore some initial ideas on how the framework we have developed in this text in the specific case of the golden ratio can be extended and generalized to other irrational bases $\gamma$ and we see that, in fact, there exist better choices for the base when constructing Hammersley sets over the golden ratio. However to finish this section, we reiterate that the convergence shown in Figure \ref{fig:star disc} (and Figure \ref{fig:L2 star disc}) when the discrepancy values are normalized by $\log N/N$ is strong evidence of the following.

\begin{conjecture}
    Given $m\geq 1$, the Hammersley point set $H_m$ in base $\varphi$ containing $N=F^m$ points, we have $$L_{2,N}(H_m) = \Theta\left(\frac{\log N}{N}\right)$$ \text{and} $$D_N^{*}(H_m) = \Theta\left(\frac{\log N}{N}\right).$$
\end{conjecture}

A significant milestone in discrepancy theory was the introduction of orthogonal function systems to analyse point sets and sequences; we refer to the excellent survey by D. Bilyk \cite{BILYK2014}. The special case of the Walsh function system is naturally used to analyse integer based digital nets and sequences due to the intrinsic link in how both objects are defined. 

\vspace{0.3cm}
\noindent
\textbf{Question 2.} Can one adapt the Walsh function system, or Haar system, to analyse and formally prove discrepancy results for the Hammersley point set in base $\varphi$?

\section{Other bases}\label{sec:otherbases}

In this final section, we discuss some initial thoughts and ideas on extending this work to constructing point sets and sequences with other irrational bases and in addition, provide the framework to study their equidistribution properties in that base. While the golden ratio without doubt has some success as a base for nets and sequences, it does however have some shortcomings. Firstly, it is not at all trivial how one should construct a $(0,2)-$sequence in base $\varphi$, or even a $(t,2)-$sequence for any parameter $t$. The second is that there is no known analogous algorithm to that of integer scrambling of nets which recovers the desirable properties of Owen scrambling. Therefore in the hope of finding a base in which the problems listed above are no more, and hopefully as a catalyst for future research, we briefly explore how to further generalize the notions of the digital point set and sequence construction and equidistribution to an irrational base $\gamma$. Throughout this last section, we will always take $\gamma$ to be the largest root of $x^2-px-q$ for $1 \le q \le p$.  The digits of numbers written in base $\gamma$ consist of $0,\dots,p$ and the expansion can always be written in \textit{reduced form} so that the digit to the right of an $p$ is strictly less that $q$. Therefore, as a general convention throughout the rest of this text, for a natural number $n$ written in base $p+1$ $$n = (d_{m-1} d_{m-2} \ldots d_1 d_0 )_{p+1}$$ we have $d_j \in \{0,1,\ldots,p\}$ and when $d_j=p$, we must have $d_{j-1}<q$. We refer to this condition as ``condition $R_{p,q}$'' and we say $n\in\mathbb{N}$ is a whole number in base $p+1$. We then use the same notation as in previous sections in this text to denote the whole number in base $\gamma$,
$$\overline{n}  = (d_{m-1}d_{m-2}\ldots d_1 d_0)_{\gamma} = \sum_{j=0}^{m-1}
d_j \gamma^j.$$

We give a general one-dimensional sequence and two-dimensional point set constructions once again following the prototypical examples of the van der Corput sequence and the Hammersley set. We study elementary intervals in base $\gamma$ and subsequently define equidistribution in base $\gamma$. We close with numerical experiments showing that the golden ratio is, in fact, not the `best' base and several open questions are posed for future research.

\subsection{One-dimensional Construction}
We form a one-dimensional sequence in base $\gamma$ by writing down the natural numbers in increasing order that satisfy the condition that any digit to the {\em left} of a $p$ in its base $p+1$ expansion is strictly less than $q$; we refer to this condition as ``condition $L_{p,q}$". (The reason why the condition is now to the left is so that when we will reverse the digits to create the van der Corput sequence in base $\gamma$, we will get an expansion that is in reduced form.) By doing so we get a sequence of natural numbers, $0=n_0<n_1<n_2<\dots$ and from here we can finish constructing the sequence in the usual manner.

\begin{defi}
Let $\gamma$ be the largest root from $x^2-px-q$ for some choice of $1 \leq q \leq p$, and let $n_i = (d_{m-1}\dots d_1 d_0)_{p+1}$ be the $i^{th}$ natural number which satisfies condition $L_{p,q}$. Then the $i^{th}$ point of the van der Corput sequence in base $\gamma$ is given by $g_{n_i}=(.d_0 d_1\dots d_{m-1})_\gamma$. 
\end{defi}

\begin{example}\label{ex:vandercorput - 1plusroot2}
    Take $p=2, q=1$ which in turn yields $\gamma = 1 + \sqrt{2}$. The following list of natural numbers satisfy condition $L_{2,1}$ and subsequently yield the first few term of the van der Corput sequence in base $1+\sqrt{2}$:
    \begin{eqnarray*}
        n_0 = 0 = (00)_{3} &\mapsto& g_0 = (0)_\gamma = 0 \\
        n_1 = 1 = (01)_{3} &\mapsto& g_1 = (.10)_\gamma = 0.4142... \\
        n_2 = 2 = (02)_{3} &\mapsto& g_2 = (.20)_\gamma = 0.8284... \\
        n_3 = 3 = (10)_{3} &\mapsto& g_3 = (.01)_\gamma = 0.1715... \\
        n_4 = 4 = (11)_{3} &\mapsto& g_4 = (.11)_\gamma = 0.5857... \\
        n_5 = 6 = (20)_{3} &\mapsto& g_5 = (.02)_\gamma = 0.3431... \\
        n_6 = 7 = (21)_{3} &\mapsto& g_6 = (.12)_\gamma = 0.7573... \\
    \end{eqnarray*}
\end{example}

\begin{figure}[t!]
    \centering
    \begin{tikzpicture}

        \draw[thick, black] (0,0) -- (10,0);

        \draw[thin,black] (0,0.2) -- (0,-0.2);
        \draw[thin,black] (10,0.2) -- (10,-0.2);
        
        \filldraw[black] (0,0) circle (3pt);
        \filldraw[black] (4.142,0) circle (3pt);
        \filldraw[black] (8.284,0) circle (3pt);
        \filldraw[black] (1.715,0) circle (3pt);
        \filldraw[black] (5.857,0) circle (3pt);
        \filldraw[black] (3.431,0) circle (3pt);
        \filldraw[black] (7.573,0) circle (3pt);

        \node at (0,-0.5) {$0$};
        \node at (10,-0.5) {$1$};
        \node at (4.142, -0.5) {$\gamma^{-1}$};
        \node at (8.284,-0.5) {$2\gamma^{-1}$};
        \node at (1.715,-0.5) {$\gamma^{-2}$};

    \draw[thick, black] (0,-1.5) -- (10,-1.5);

        \draw[thin,black] (0,-1.3) -- (0,-1.7);
        \draw[thin,black] (10,-1.3) -- (10,-1.7);
        
        \filldraw[black] (0,-1.5) circle (3pt);
        \filldraw[black] (4.142,-1.5) circle (3pt);
        \filldraw[black] (8.284,-1.5) circle (3pt);
        \filldraw[black] (1.715,-1.5) circle (3pt);
        \filldraw[black] (5.857,-1.5) circle (3pt);
        \filldraw[black] (3.431,-1.5) circle (3pt);
        \filldraw[black] (7.573,-1.5) circle (3pt);

        \filldraw[black] (0,-1.5) circle (3pt);
        \filldraw[black] (0.747,-1.5) circle (3pt);
        \filldraw[black] (4.852,-1.5) circle (3pt);
        \filldraw[black] (8.994,-1.5) circle (3pt);
        \filldraw[black] (2.426,-1.5) circle (3pt);
        \filldraw[black] (6.568,-1.5) circle (3pt);
        \filldraw[black] (1.421,-1.5) circle (3pt);
        \filldraw[black] (5.563,-1.5) circle (3pt);
        \filldraw[black] (9.705,-1.5) circle (3pt);
        \filldraw[black] (3.137,-1.5) circle (3pt);
        \filldraw[black] (7.279,-1.5) circle (3pt);

    \end{tikzpicture}
    \caption{The first 7 and 17 terms in the van der Corput sequence in base $\gamma = 1+\sqrt{2}$.}
    \label{fig:vdcotherbaseexample - first 7 terms}
\end{figure}
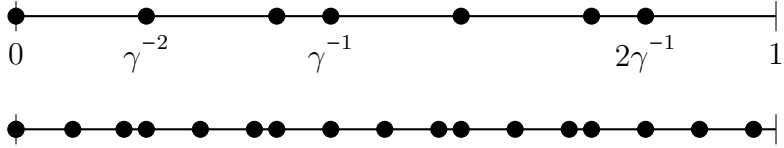

Instead of relying on an ad-hoc method of constructing the sequence $0=n_0<n_1<n_2<\dots$ of natural numbers satisfying condition $L_{p,q}$, we provide the following recursive algorithm.

\begin{algorithm}[Condition $L_{p,q}$]

To generate the set of $n\in \mathbb{N}$ satisfying condition $L_{p,q}$:
\begin{enumerate}
     \item Use the convention that there is exactly one number (namely $0$) with zero digits.
     \item The numbers whose base $p+1$ expansion has one digit and satisfies the condition $L_{p,q}$ are $0,1,\dots,p$. There are $p+1$ of these.
     \item Let $G_{m-1}$ be the number of $(m-1)-$digit numbers satisfying the condition $L_{p,q}$, and let $A_{m-1}$ be the number of $(m-1)-$digit numbers that satisfy the condition $L_{p,q}$ and whose most significant digit is not $p$. Finally, let $\Gamma^L_{m}=\{0=n_0<\dots<n_{G_m - 1}\}$ be the natural numbers satisfying condition $L_{p,q}$ with at most $m$ digits in increasing order. The first $G_{m-1}$ numbers in $\Gamma^L_{m}$ are just those from $\Gamma^L_{m-1}$ (i.e., those whose most significant digit is $0$). The numbers $n_i$ in $\Gamma^L_{m}$ whose $m$th digit is $d\in\{1,\dots,q-1\}$ are those with indices $i$ between $dG_{m-1}$ and $(d+1)G_{m-1}-1$ (inclusive) and of the form
    \[
     n_i=d(p+1)^{m-1}+n_{j}
      \]
      with $i=dG_{m-1}+j$ where $0\leq j<G_{m-1}$. The numbers in $\Gamma^L_{m}$ whose most significant digit is $d\in\{q,q+1,\dots,p\}$ are those with indices $i$ between $qG_{m-1}+(d-q)A_{m-1}$ and $qG_{m-1}+(d-q+1)A_{m-1}-1$ (inclusive) with 
      \[
      n_i=d(p+1)^{m-1} + n_j,
      \]
      and $i=qG_{m-1}+j$ where $0\leq j< A_{m-1}$. Thus we have the recurrences for the cardinality of the set $\Gamma_m^{L}$:
     \begin{align*}
          &G_0 = 1,\quad  G_1 = p+1,  &A_0 = 1,\quad A_1 = p,\\
         &G_m = qG_{m-1} + (p-q+1)A_{m-1},&A_{m}= qG_{m-1} + (p-q)A_{m-1}.
     \end{align*}
 \end{enumerate}
\end{algorithm}

\begin{example}
Extending from Example \ref{ex:vandercorput - 1plusroot2}, using the above algorithm we see that for base $1+\sqrt{2}$ the numbers $G_m$ are as follows: $G_0 = 1, G_1 = 3, G_2 = 7, G_3 = 17, G_4 = 41, G_5 = 99, ...$
\end{example}

\subsection{Two-dimensional construction}


We now want to construct a finite point set in two dimensions and the usual Hammersley construction seems like a reasonable place to begin. Keeping the one-dimensional van der Corput sequence in base $\gamma$ as defined in the previous section as the first coordinate, we want to create a two-dimensional point set with $G_m$ points by adding a second coordinate. To do so, we firstly gather the $m-$digit whole numbers in base $\gamma$. That is, we find those natural numbers which have a finite base $p+1$ expansion satisfying condition $R_{p,q}$ and denote these by $\Gamma^R_{m} = \{ 0 = n_0 < n_1 < n_2 < \cdots \}$. As remarked much earlier in this text, one way to describe the $i^{th}$ point of the traditional Hammersley construction is to define the first coordinate as the $i^{th}$ term of the van der Corput sequence, and then the second coordinate is obtained by reversing the order of the digits of the number in the first coordinate. In general, in base $\gamma$ this exact construction is not possible. Despite the fact that the cardinalities of the sets $\Gamma^L_{m}$ and $\Gamma^R_{m}$ are equivalent, the natural numbers contained in the sets are not and therefore the digit expansions are not "reversible". It can be shown that the sets $\Gamma^L_{m}$ and $\Gamma^R_{m}$ are equivalent if and only if $p=q$.

\begin{figure}[t!]
    \centering
    \includegraphics[width=7cm]{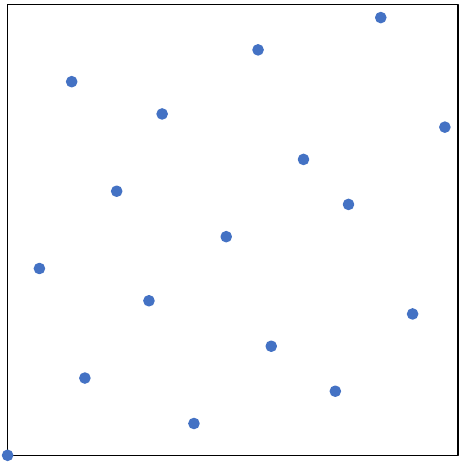}
    \caption{The 3-digit Hammersley point set in base $1+\sqrt{2}$.}\label{fig:hammersley 1 plus root 2 - 17 points}
\end{figure}

Before formally defining the two-dimensional construction, we give a second algorithm using a recursive approach similar to the one used in the previous section to find the sequence of natural numbers contained in $\Gamma^R_{m}$.

\begin{algorithm}[Condition $R_{p,q}$]\label{algo:right}
To generate the set of $n\in\mathbb{N}_0$ satisfying condition $R_{p,q}$:
 \begin{enumerate}
     \item The only $0-$digit number whose base $p+1$ expansion satisfies the condition $R_{p,q}$ is $0$.
     \item The $1-$digit numbers whose base $p+1$ expansion satisfies the condition $R_{p,q}$ are $0,1,\dots,p$.
     \item Observe that reversing the digits of the base $p+1$ expansion of the numbers whose base $p+1$ expansion has at most $m$ digits and satisfy $L_{p,q}$ is a re-ordering of the base $p+1$ expansions of the numbers whose base $p+1$ expansion has at most $m$ digits and satisfy the condition $R_{p,q}$. Thus $G_m$ counts the number of numbers in base $p+1$ with at most $m$ digits and satisfying $R_{p,q}$. Let $\Gamma^R_{m}$ be the natural numbers (in increasing order) whose base $p+1$ expansions have at most $m-$digits and satisfy the condition $R_{p,q}$. Let $B_m$ count the number of numbers in $\Gamma^R_{m}$ whose $m^{th}$ digit is strictly less than $q$.
    \item The first $G_{m-1}$ numbers in $\Gamma^R_{m}$ are simply those from $\Gamma^R_{m-1}$ (their $m^{th}$ digit is $0$). The numbers in $\Gamma^R_{m}$ whose $m^{th}$ digit is $d\in\{1,\dots,p-1\}$ are those with indices $dG_{m-1}$ to $(d+1)G_{m-1}-1$ (inclusive), and the $i^{th}$ number with $i=dG_{m-1}+j$, where $0\leq j < G_{m-1}$ is 
    \[
     n_i=n_j+d(p+1)^{m-1}.
     \]
     The numbers in $\Gamma^R_{m}$ whose $m^{th}$ digit is $p$ are those with indices $pG_{m-1}\leq i<pG_{m-1}+B_{m-1}$ and 
     \[
     n_i=n_j+p(p+1)^{m-1}.
     \]
     We get the recurrences:
         \begin{align*}
          &G_0 = 1,\quad  G_1 = p+1,  &B_0 = 1,\quad B_1 = q,\\
        &G_m = pG_{m-1} + B_{m-1},&B_{m}= qG_{m-1}
     \end{align*}
     Thus $G_m=pG_{m-1}+qG_{m-2}$.
 \end{enumerate}
\end{algorithm}

We can now define the two-dimensional Hammersley point set in base $\gamma$.

\begin{defi}\label{def:basegammatwodimensional}
Let $\gamma$ be the largest root from $x^2-px-q$ for some choice of $1 \leq q \leq p$. Then the $m-$digit Hammersley point set $P_m$ in base $\gamma$ containing $G_m$ points is defined as $$H^\gamma_m := \left\{\left( g_{n_i},\frac{\overline{n_i}}{\gamma^m} \right) : 0 \leq i < G_m \right\}$$ where $g_{n_i}$ is the $i^{th}$ term of the van der Corput sequence in base $\gamma$ and $n_i \in \Gamma^R_{m}$.
\end{defi}

\subsection{Equidistribution in other bases}\label{sec:equidist in other bases}

We can begin by defining the $m-$partition in base $\gamma$. Given $\Gamma^R_{m} = \{0=n_0<n_1<\dots<n_{G_m-1}\}$ an elementary interval in base $\gamma$ is of the form
 \[
  \Bigg[\frac{\overline{n_l}}{\gamma^m},\frac{\overline{n_{l+1}}}{\gamma^m}\Bigg),
 \]
where $n_l \in \Gamma_m^R$. For a fixed $m$, the set of all such intervals forms a partition of $[0,1)$ which we call the $m-$partition in base $\gamma$ and denote this object $\mathcal{P}^{\gamma}_{m}$.

\begin{defi}
    An \textit{elementary $m-$interval in base $\gamma$} is a subinterval of $[0,1)$ of the form $$I = \Bigg[\frac{\overline{n_l}}{\gamma^m},\frac{\overline{n_{l+1}}}{\gamma^m}\Bigg)$$ for some $n_l \in \Gamma^R_{m}$.
\end{defi}

We observe, as in the case of the golden ratio base, that elementary intervals of the same size appear in several different partitions of the unit interval. For this reason, we must distinguish between elementary intervals in base $\gamma$ and \textit{prime} elementary intervals in base $\gamma$. This is fulfilled by Notation \ref{not:I notation other bases} in an analogous manner to the golden ratio base.

In addition to characterizing the prime elementary intervals, we need to determine the size of the elementary $m-$intervals in base $\gamma$ and subsequently the number of points from 
a $G^m$ size point set $P_m$ that should be required in each elementary interval. From here, one can define the conditions and properties required to study the equidistribution of $P_m$ in base $\gamma$. The following lemma takes the first of these steps giving the possible lengths of elementary intervals in base $\gamma$.

 \begin{lemma}\label{lem:length of base gamma intervals}
     For $n_l \in \Gamma^R_{m}$, an elementary $m-$interval $I = \left[\frac{\overline{n_l}}{\gamma^m},\frac{\overline{n_{l+1}}}{\gamma^m}\right)$ in base $\gamma$ has two possible lengths of $\gamma^{-m}$ or $q\gamma^{-m-1}$ depending on the last digit of the base $\gamma$ expansion of $n_l$. In particular, if $|n_l|<p$ then $I$ has length $\gamma^{-m}$ and if $|n_l| = p$ then $I$ has length $q\gamma^{-m-1}$.
 \end{lemma}

 \begin{proof}
     There are two cases to consider. The first case occurs when the smallest digit of $n_l$ is strictly less than $p$, i.e. suppose that $n_l=(d_{m-1}\cdots d_1d_0)_{p+1}$ with $0\leq d_0<p$. Then $n_{l+1}=(d_{m-1}\cdots d_1(d_0+1))_{p+1}$ so that the interval is of the form
 \[
 \left[\frac{\overline{n_l}}{\gamma^m},\frac{\overline{n_{l+1}}}{\gamma^m}\right) = \left[\frac{\overline{n_l}}{\gamma^m}, \frac{\overline{n_l}}{\gamma^m}+\frac{1}{\gamma^m}\right)
 \]
 and the length of the interval is clearly $\gamma^{-m}$. For the second case, let the smallest digit of $n_l$ equal $p$. To see the result, we use the relation from the minimal polynomial $x^2-px-q$ of $\gamma$ to get
 \[
 \overline{n_l}=(d_{m-1}\cdots d_2d_1p)_{\gamma} \text{ and } \overline{n_{l+1}}=(d_{m-1}\cdots d_2(d_1+1)0)_{\gamma} = (d_{m-1}\cdots d_2d_1p)_{\gamma} + q\gamma^{-1}
 \]
 so that the interval is of the form
 \[
 \left[\frac{\overline{n_l}}{\gamma^m},\frac{\overline{n_l}}{\gamma^m}+\frac{q}{\gamma^{m+1}}\right),
 \]
 and it follows that the length of the interval is $q\gamma^{-m-1}$.
 \end{proof}

\begin{remark}
We note that in the case of $p=q=1$ (when $\gamma = \varphi$), this notion agrees with the content earlier in the paper when elementary intervals in base $\varphi$ had the permitted lengths of $\varphi^{-m}$ and $(1)\varphi^{-m-1}=\varphi^{-m-1}$. After performing extensive numerical experiments, it appears that the choice $q=1$ yields the best distribution properties of the resulting two-dimensional point sets; we encourage the reader to study Tables \ref{table:p2q1first} to \ref{tab:p4q4last} in the appendix for empirical evidence of this claim. Thus, for the remainder of this text we assume $q=1$.  Consequently, from Lemma \ref{lem:length of base gamma intervals}, the only possible lengths of elementary intervals in base $\gamma$ are $\gamma^{-m}$ and $\gamma^{-m-1}$. 
\end{remark}

We can introduce analogous notation as in Notation \ref{not:elementary intervals} for the more general base $\gamma$ allowing us to fulfil the important tasks of formalising the orientation of the various lengths of elementary intervals and thus characterizing all prime elementary intervals as discussed earlier in the paper.

\begin{notation}\label{not:I notation other bases}
    Assume $p \ge 2$. Let $I$ be an elementary interval in the $m-$partition in base $\gamma$ with $$I = \left[\frac{\overline{n_l}}{\gamma^m}, \frac{\overline{n_{l+1}}}{\gamma^m}\right),$$ where $n_l \in \Gamma^R_{m}$. Set $n_a = (d_{m-1}\ldots d_2 d_1 0)_{p+1}$ and choose $i\in\{0,1, \ldots, p+1\}$ such that $n_l = n_{a+i}$. We use the notation $I^\gamma_i(a;m)$ to denote $I$, i.e. $$I = I_i^\gamma(a;m) \coloneqq \left[ \frac{\overline{n_{a+i}}}{\gamma^m}, \frac{\overline{n_{a+i+1}}}{\gamma^m} \right).$$
\end{notation}

Explicitly, $n_a = (d_{m-1}\ldots d_2 d_1d_0)_{p+1}$ with $d_1 \in \{0,1,\ldots, p-2\}$ implies that only the intervals $I^\gamma_0(a;m), \ldots, I^\gamma_{p}(a;m)$ exist and $n_a = (d_{m-1}\ldots d_2 (p-1)0)_{p+1}$ implies that all intervals $I^\gamma_0(a;m), \ldots, I^\gamma_{p+1}(a;m)$ exist. We refer the reader back to Figure \ref{fig:vdcotherbaseexample - first 7 terms} to study the $m-$partition in base $\gamma = 1 + \sqrt{2}$ ($p=2,q=1$) for $m=2$ and $m=3$. 

\begin{remark}
There are $G_m$ elementary intervals in $\mathcal{P}_m^\gamma$. This fact follows directly from noting the cardinality of the set $\Gamma_m^R$ as established in Algorithm \ref{algo:right}. An elementary interval $$I = \left[ \frac{\overline{n_l}}{\gamma^m}, \frac{\overline{n_{l+1}}}{\gamma^m} \right)$$ with $n_l \in \Gamma^R_{m}$ is a prime elementary interval if and only if $d_1 \neq p$ (where here $d_1$ denotes the second-to-last digit of $n_l$), i.e. if and only if $I = I^\gamma_i(a;m)$ for $i \in \{0,1,\ldots,p\}$.\hnathan{prime $\iff d_1 \neq p$ ??}
\end{remark}

\begin{defi}
    A \textit{prime elementary $(k_1,\ldots,k_s)-$interval} in base $\gamma$ is a subset of $[0,1)^s$ of the form $$I = \prod_{j=1}^s \left[ \frac{\overline{n_{l,j}}}{\gamma^{k_j}}, \frac{\overline{n_{l+1,j}}}{\gamma^{k_j}} \right)$$ where $n_{l,j} \in \Gamma_{k_j}^R$ has its second-to-last digit not equal to $p$. 
\end{defi}

\begin{notation}
    For a prime elementary $\mathbf{k} = (k_1,\ldots,k_s)$-interval $I$ based on a vector $\mathbf{k} \in \mathbb{N}_0^s$, we use the notation $$|I| = \sum_{j=1}^s k_j + l_p$$ where $l_p$ counts the number of $1 \leq j \leq s$ such that $|n_{l,j}|=p$. This convention implies that the volume of $I$ is $\gamma^{-|I|}.$
\end{notation}

The following series of lemmas allow us to deduce how many points from a total of $G_m$ should be contained within an elementary interval in base $\gamma$ of a given volume. The analogous result in the golden ratio case is Lemma \ref{lem:bestapproximation}.

\begin{lemma}\label{lem:closedgm}
    For $m\geq0$, the numbers $G_m$ from Algorithm \ref{algo:right} have a closed form expression given as 
    \begin{equation}\label{eq:gm}
    G_m = \frac{(\gamma+1) \gamma^m + (\gamma - p - 1)(p-\gamma)^m}{2\gamma-p}
    \end{equation}
    where $\gamma$ is the largest root of $x^2-px-q$ for $1\leq q\leq p$.
\end{lemma}
\begin{proof}
    Beginning with the relation $G_m = pG_{m-1}+qG_{m-2}$ as derived in Algorithm \ref{algo:right}, we arrive at the recurrence characteristic polynomial $x^2-px-q$. We know that one root (the largest of the two) of this polynomial is $\gamma$ and the other, by standard quadratic polynomial algebra, can be expressed as $p - \gamma$. This gives a general solution of $$G_m = c_1 \gamma^m + c_2 (p-\gamma)^m$$ where the constants $c_1$ and $c_2$ are to be determined by initial terms of the recurrence. In our case, $G_0 = 1$ and $G_1 = p+1$ producing 
    \begin{eqnarray*}
        1 &=& c_1+c_2 \\ 
        p+1 &=& c_1 \gamma + c_2 (p-\gamma).
    \end{eqnarray*}
    Solving this system for $c_1$ and $c_2$ yields the formula in the statement of the lemma.
\end{proof}
\begin{remark}
    For completeness, with $p = q = 1$ and subsequently $\gamma = \varphi$, the formula \eqref{eq:gm} in turn reproduces \eqref{eq: fib as phi} from Section \ref{sec:definitionsandnotation}. Explicitly, for $m\geq0$ we have $G_m = F^m = \frac{\varphi^{m+2}-\psi^{m+2}}{\sqrt{5}}$ for $\gamma = \varphi$.
\end{remark}
Next, assuming that we are given a point set which contains $G_m$ total points for $m\geq0$, we are confronted with the same problem as before created by the volume of an elementary interval in base $\gamma$ being irrational. We shall never be able to force the (local) discrepancy of an elementary interval to be zero and the best we can hope to do is a minimization.
\begin{lemma}
    For $m, k\geq 0$ with $m \geq k$, $$\min_{0\leq i\leq m} \left| \frac{G_{m-i}}{G_m} - \frac{1}{\gamma^k} \right| = \left| \frac{G_{m-k}}{G_m} - \frac{1}{\gamma^k} \right|.$$
\end{lemma}
\begin{proof}
    From the expression derived in Lemma \ref{lem:closedgm}, we write 
    \begin{eqnarray*}
        \left| \frac{G_{m-i}}{G_m} - \frac{1}{\gamma^k} \right| &=& \left| \frac{(\gamma+1) \gamma^{m-i} + (\gamma - p - 1)(p-\gamma)^{m-i}}{(\gamma+1) \gamma^{m} + (\gamma - p - 1)(p-\gamma)^{m}} - \frac{1}{\gamma^k} \right| \nonumber \\
        &=&\frac{1}{G_m} \left| (\gamma-p-1)((p-\gamma)^{m-k}-(p-\gamma)^m)\right|\nonumber \\
        &=&\left|\frac{(p-\gamma)^{m-k}}{G_m}\right| \left| (\gamma-p-1)) (1-\gamma^{-k}(-1)^k)\right|\nonumber \\      
    \end{eqnarray*}
    We note that $\left|(p-\gamma)^{m-k}\right|<1$. Furthermore, $(1-\gamma^{-k}(-1)^k)>0$ is clearly maximized when $k=1$, in which case it equals $1+1/\gamma=1+\gamma-p$. Hence
    \begin{eqnarray*}
        \left| \frac{G_{m-i}}{G_m} - \frac{1}{\gamma^k} \right|
        &<& \frac{1}{G_m} \left|(\gamma-p-1)(\gamma-p+1) \right| = \frac{1}{G_m} \left|(\gamma-p)^2-1)\right| \\
        &=& \frac{1}{G_m} \left|p(p-\gamma)\right|= \frac{1}{G_m} \left|\frac{p}{\gamma}\right| <\frac{1}{G_m}.
    \end{eqnarray*}  
    Hence $G_{m-k}/G_m$ is the best rational approximation with denominator $G_m$ for $\gamma^{-k}$.
\end{proof}

From the above exposition, we can conclude that an elementary interval of size $\gamma^{-k}$ should contain $G_{m-k}$ points from a point set of size $G_m$ in order for $P_m$ to minimize the local discrepancy and achieve good distribution properties. More formally, we define equidistribution, nets and sequences in base $\gamma$.

\begin{figure}[t!]
    \centering
    \includegraphics[width=15cm]{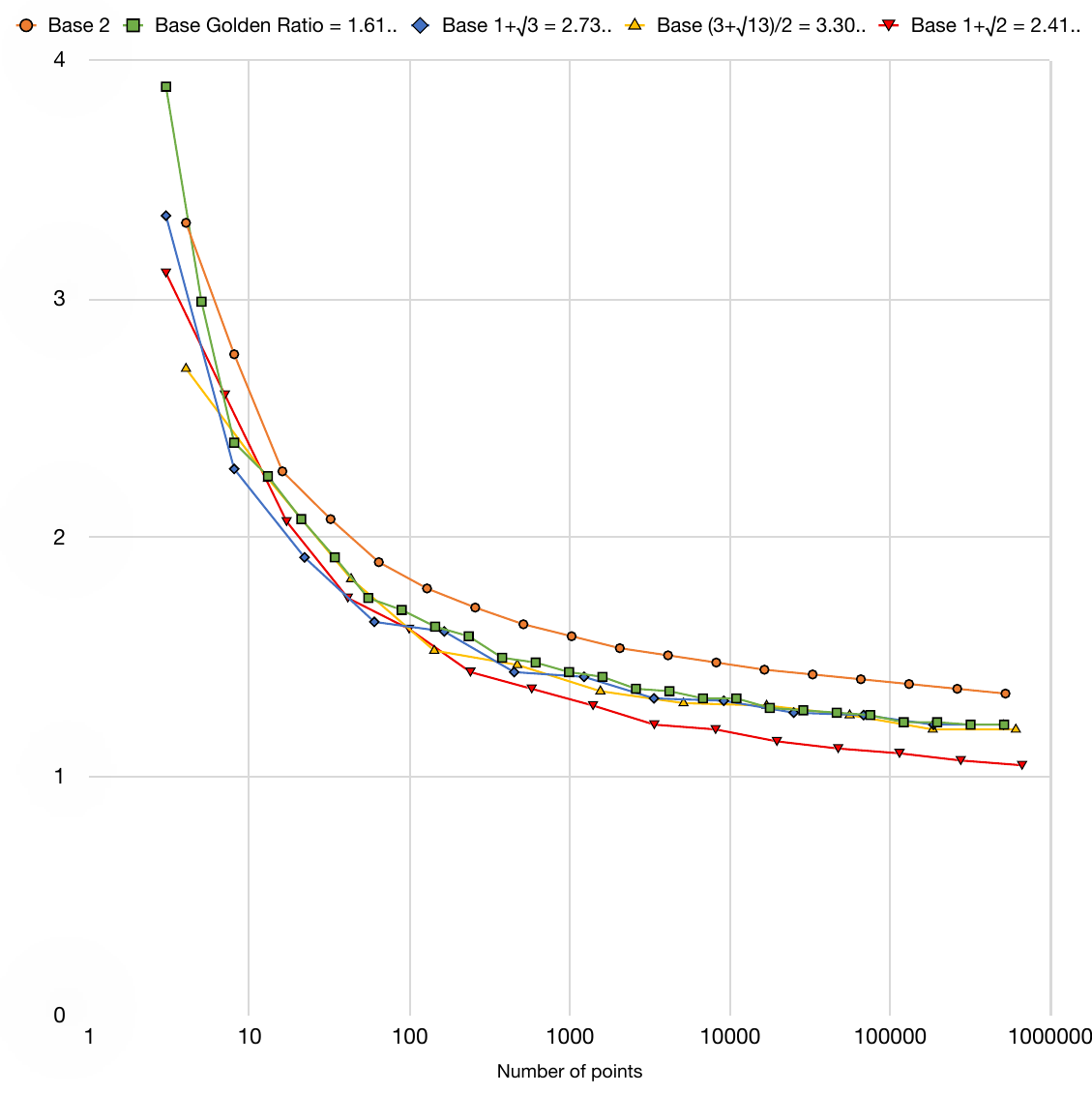}
    \caption{The best performing two-dimensional Hammersley point sets in base $\gamma$ with respect to the star discrepancy. All values are normalized by $\log N/N$.}
    \label{fig:discrepancy plot other bases}
\end{figure}

\begin{defi}
    A point set $P_m$ contained in $[0,1)^s$ with $G_m$ points is said to be \textit{$(k_1,\ldots,k_s)$-eq\-ui\-dis\-tri\-buted in base $\gamma$} if every prime elementary $(k_1,\ldots,k_s)-$interval in base $\gamma$, $I$, contains exactly $G_{m-|I|}$ points from $P_m$. 
\end{defi}

Since every elementary interval in $\mathcal{P}_{\mathbf{m}}^\gamma$ is not classified as a prime elementary interval, we have the following stronger form of equidistribution.

\begin{defi}
    A point set $P_m$ contained in $[0,1)^s$ with $G_m$ points is said to be \textit{strongly $(m_1,\ldots,m_s)-$equidistributed in base $\gamma$} if every elementary interval $I$ in the $(m_1,\ldots,m_s)$-partition in base $\gamma$ contains exactly $G_{m-|I|}$ points from $P_m$. 
\end{defi}

\begin{defi}
    A point set $P_m$ contained in $[0,1)^s$ with $G_m$ points is called a $(t,m,s)-$net in base $\gamma$ if it is $(k_1,\ldots,k_s)-$equidistributed in base $\gamma$ for all $\rho(\mathbf{k}) \leq m - t$.
\end{defi}

To finish, we present some extended numerical results in Figure \ref{fig:discrepancy plot other bases} for the star discrepancy of Hammersley point set in base $\gamma$ as defined in Definition \ref{def:basegammatwodimensional} for several ``good" choices of $\gamma.$

\subsection*{Future Work} The authors hope that this paper acts as a catalyst for further investigation into utilizing irrational bases to construct digital nets and sequences. We now gather and present a list of open questions for future research:
\begin{enumerate}
    \item Under what base $\gamma \in \mathbb{R}$, does the Hammersley construction as presented in Definition \ref{def:basegammatwodimensional} achieve the lowest star discrepancy? Empirically, the ``best" base found so far is $1+\sqrt{2}$ as per Figure \ref{fig:discrepancy plot other bases}.
    \item Based on numerical experiments of the discrepancy of irrational based digital nets and sequences, why is a base chosen between $2$ and $3$ seemingly the best base for construction?
    \item In the context of choosing the largest root from the polynomial $x^2-px-q$, are the best irrationals for construction those that are a bit bigger than $p$, a little less than $p+1$, or halfway between $p$ and $p+1$?
    \item Is it possible that an adaptation can be made to the Haar or Walsh function systems to employ an orthogonal function method to derive formal bounds for the discrepancy of the van der Corput sequence or Hammersley point set in base $\varphi$?
    \item Can one construct a scrambling algorithm in base $\varphi$ which preserves the desirable properties of Owen scrambling? Alternatively, from Conjecture \ref{con:quasi-equidistribution}, numerical exploration suggests that one can scramble base $\varphi$ point sets in base 2 and preserve good distribution properties. Hence, how do the RQMC estimators derived from the base 2 scrambled golden ratio Hammersley net compare to the traditional base 2 Hammersley estimator?
    \item How should one construct digital nets and sequences in base $\varphi$ or general base $\gamma$ in higher than dimension two?
    \item Can one derive existence results for irrational based digital nets and sequences? Specifically, can one find a $(t,2)-$sequence construction in base $\varphi$ for any $t \geq 0$?

\end{enumerate}

\section*{Acknowledgements} The authors would like to thank Fran\c{c}ois Cl\'{e}ment for providing the DEM algorithm (with original code by Markus Wahlstrom) allowing faster computation of the star discrepancy in our numerical calculations.

\bibliographystyle{plain}
\bibliography{refs}

\begin{bibdiv}
\begin{biblist}

\bib{BILYK2014}{incollection}{
      author={Bilyk, D.},
       title={Roth's orthogonal function method in discrepancy theory and some new connections},
        date={2014},
   booktitle={A panorama of discrepancy theory},
      series={Lecture Notes in Math.},
      volume={2107},
   publisher={Springer, Cham},
       pages={71\ndash 158},
}

\bib{CARBONE2012}{article}{
      author={Carbone, I.},
       title={Discrepancy of {LS}-sequences of partitions and points},
        date={2012},
     journal={Ann. Mat.},
      volume={191},
       pages={819\ndash 844},
}

\bib{CARBONE2016}{inproceedings}{
      author={Carbone, I.},
       title={Comparison between {LS}-sequences and $\beta-$adic van der {C}orput sequences},
        date={2015},
   booktitle={Monte {C}arlo and quasi-{M}onte {C}arlo methods in scientific computing},
   publisher={Springer},
       pages={261\ndash 270},
}

\bib{CARBONEIACOVOLCIC2012}{article}{
      author={Carbone, I.},
      author={Iac\'o, M.~R.},
      author={Vol\v{c}i\v{c}, A.},
       title={{LS}-sequences of points in the unit square},
        date={2012},
         url={https://arxiv.org/abs/1211.2941},
}

\bib{DICKPILL2010}{book}{
      author={Dick, J.},
      author={Pillichshammer, F.},
       title={Digital nets and sequences},
   publisher={Cambridge University Press, Cambridge},
        date={2010},
}

\bib{DRMOTATICHY1997}{book}{
      editor={Drmota, M.},
      editor={Tichy, R.F.},
       title={Sequences, discrepancies and applications},
      series={Lecture Notes in Mathematics},
   publisher={Springer-Verlag},
        date={1997},
      volume={1651},
}

\bib{FAURE1981}{article}{
      author={Faure, H.},
       title={Discr\'{e}pance de suites associ\'{e}es \`{a} un syst\`eme de num\'{e}ration (en dimension un)},
    language={French},
        date={1981},
     journal={Bulletin de la Soci\'{e}t\'{e} Math\'{e}matique de France},
      volume={109},
       pages={143\ndash 182},
         url={http://www.numdam.org/articles/10.24033/bsmf.1935/},
}

\bib{FAURE1982}{article}{
      author={Faure, H.},
       title={Discr\'{e}pance de suites associ\'{e}s \`{a} un syst\'{e}me de num\'{e}ration (en dimension s)},
    language={French},
        date={1982},
     journal={Acta Arithmetica},
      volume={4},
       pages={337\ndash 351},
         url={http://eudml.org/doc/205851},
}

\bib{HLAWKA1961}{article}{
      author={Hlawka, E.},
       title={Funktionen von beschr\"{a}nkter variation in der theorie der gleichverteilung},
        date={1987},
     journal={Ann. Mat. Pura Appl.},
      volume={54},
       pages={325\ndash 333},
}

\bib{KAKUTANI1976}{inproceedings}{
      author={Kakutani, S.},
       title={A problem of equidistribution on the unit interval [0,1]},
        date={(1976)},
   booktitle={Measure theory},
      editor={Bellow, Alexandra},
      editor={K{\"o}lzow, Dietrich},
   publisher={Springer Berlin Heidelberg},
     address={Berlin, Heidelberg},
       pages={369\ndash 375},
}

\bib{KOKSMA1942}{article}{
      author={Koksma, J.~F.},
       title={A general theorem from the theory of uniform distribution modulo 1},
        date={1942},
     journal={Mathematica, Zutphen. B.},
      volume={11},
       pages={7\ndash 11},
}

\bib{LEMIEUX2009}{book}{
      author={Lemieux, C.},
       title={Monte {C}arlo and quasi-{M}onte {C}arlo sampling},
   publisher={Springer Science \& Business Media, Springer},
        date={2009},
}

\bib{LEOPILL2014}{book}{
      author={Leobacher, G.},
      author={Pillichshammer, F.},
       title={Introduction to quasi-{M}onte {C}arlo integration and applications},
   publisher={Birkh\"{a}user},
        date={2014},
}

\bib{NIE1992}{book}{
      author={Niederreiter, H.},
       title={Random number generation and quasi-{M}onte {C}arlo methods},
   publisher={SIAM, Philadelphia},
        date={1992},
}

\bib{Nin98}{article}{
      author={Ninomiya, S.},
       title={Constructing a new class of low-discrepancy sequences by using the $\beta$-adic transformation},
        date={1998},
     journal={Mathematics and Computers in Simulation},
      volume={47},
       pages={403\ndash 418},
}

\bib{NINOMIYA1998}{article}{
      author={Ninomiya, S.},
       title={On the discrepancy of the $\beta-$adic van der corput sequence},
        date={1998},
     journal={J. Math. Sci. Univ. Tokyo},
      volume={5},
       pages={345\ndash 366},
}

\bib{OWEN1995}{inproceedings}{
      author={Owen, A.~B.},
       title={Randomly permuted (t,m,s)-nets and (t,s)-sequences},
        date={1995},
   booktitle={Monte {C}arlo and quasi-{M}onte {C}arlo methods in scientific computing},
      editor={Niederreiter, Harald},
      editor={Shiue, Peter Jau-Shyong},
   publisher={Springer New York},
     address={New York, NY},
       pages={299\ndash 317},
}

\bib{NIED2014}{inproceedings}{
      author={Schretter, C.},
      author={Zhijian, H.},
      author={Gerber, M.},
      author={Chopin, N.},
      author={Niederreiter, H.},
       title={Van der {C}orput and golden ratio sequences along the hilbert space-filling curve},
        date={(2015)},
   booktitle={Monte {C}arlo and quasi-{M}onte {C}arlo methods in scientific computing},
   publisher={Springer},
       pages={531\ndash 544},
}

\bib{rSOB67a}{article}{
      author={Sobol', I.~M.},
       title={On the distribution of points in a cube and the approximate evaluation of integrals},
        date={1967},
     journal={USSR Comp. Math. Math. Phys.},
      volume={7},
       pages={86\ndash 112},
}

\bib{VOLCIC2011}{article}{
      author={Vol\v{c}i\v{c}, A.},
       title={A generalization of {K}akutani's splitting procedure},
        date={2011},
     journal={Ann. Mat.},
      volume={190},
       pages={45\ndash 54},
}

\bib{CARBONEVOLCIC2007}{article}{
      author={Vol\v{c}i\v{c}, A.},
      author={Carbone, I.},
       title={Kakutani's splitting procedure in higher dimension},
        date={2007},
     journal={Rend. Ist. Matem. Univ. Trieste},
      volume={39},
       pages={1\ndash 8},
}

\bib{WIARTLEMIEUXDONG2021}{article}{
      author={Wiart, J.},
      author={Lemieux, C.},
      author={Dong, G.~Y.},
       title={On the dependence structure and quality of scrambled (t,m,s)-nets},
        date={2021},
     journal={Monte Carlo Methods and Applications},
      volume={27},
      number={1},
       pages={1\ndash 26},
}

\end{biblist}
\end{bibdiv}

\newpage
\appendix

\section{Additional Figures - Golden Ratio Nets and Sequences}\label{app:additional figures 1}

For the benefit of the reader, further graphics of nets and sequences in golden ratio base are given below.

\begin{figure}[h!]
     \centering
     \includegraphics[scale =.20]{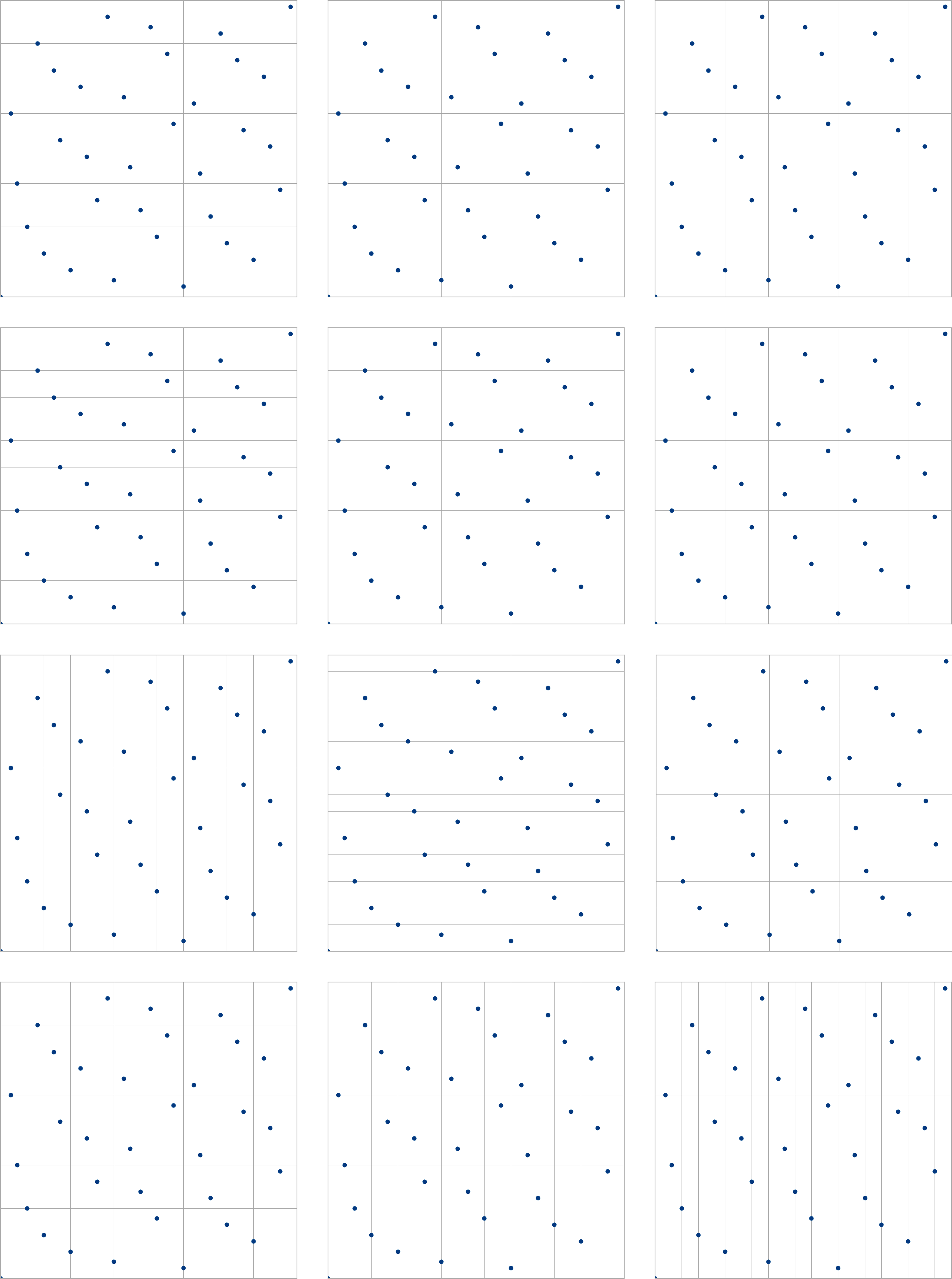}
     \caption{The Hammersley set $H_7$ shown to be a $(0,7,2)-$net with all valid $(m_1,m_2)-$partitions.}\label{fig: various partitions}
\end{figure}

\begin{figure}[h!]
     \centering
     \includegraphics[scale = .34]{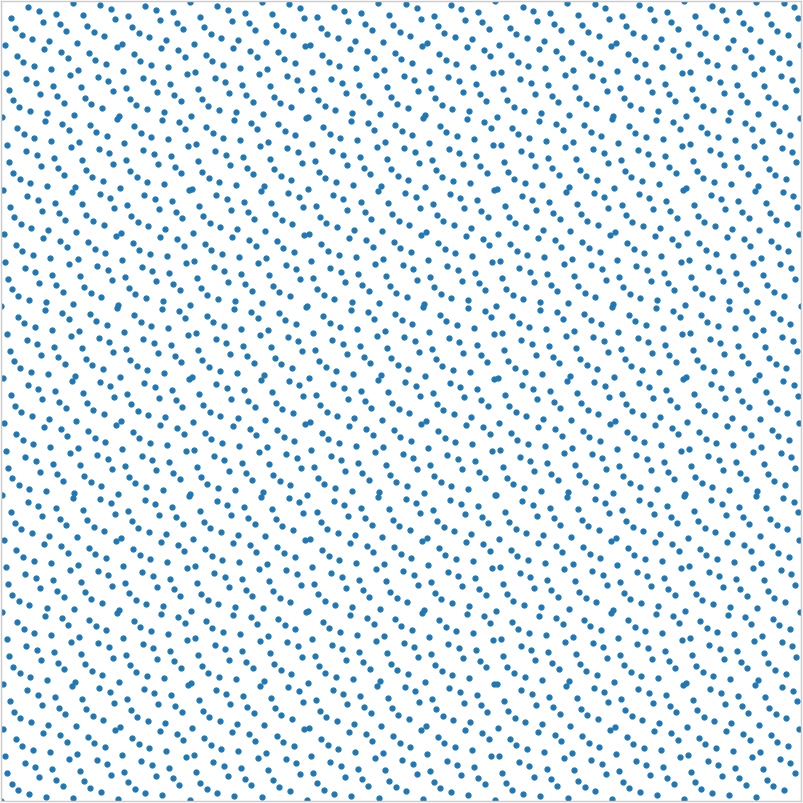}
     \caption{The Hammersley point set $H_{16}$ in base $\varphi$.}\label{fig: hammersley 16}
\end{figure}

\begin{figure}[h!]
    \centering
    \includegraphics[scale=0.15]{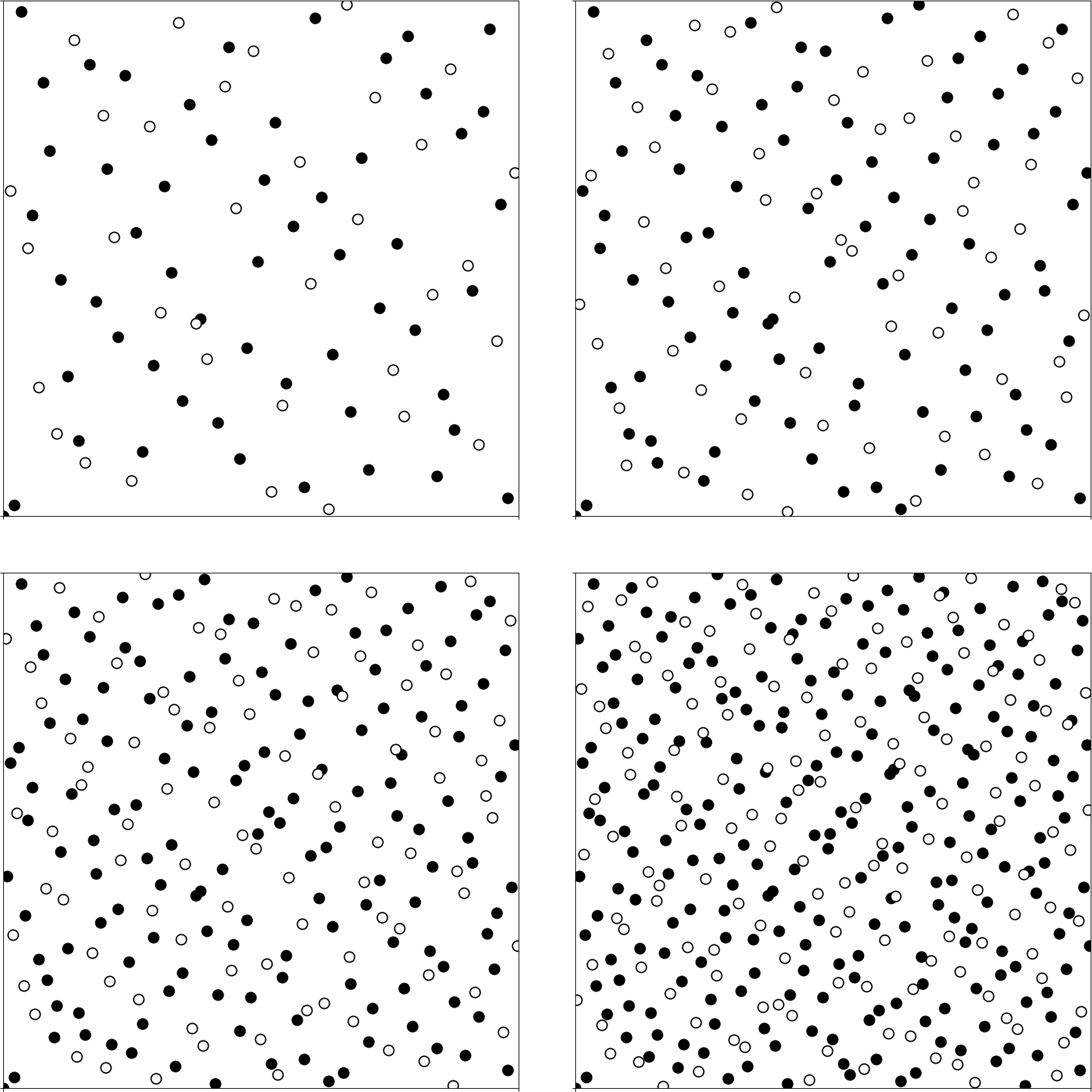}
    \caption{The sequence generated by applying Lemma \ref{lem: Extending points} repeatedly with the origin as the starting point. The points to be added at each step are circles. Given are the steps 9 to 12, extended from Figure \ref{fig:extensions from lemma steps 1 to 8} in the main text.}
    \label{fig:sequence lemma extension 9 to 12}
\end{figure}

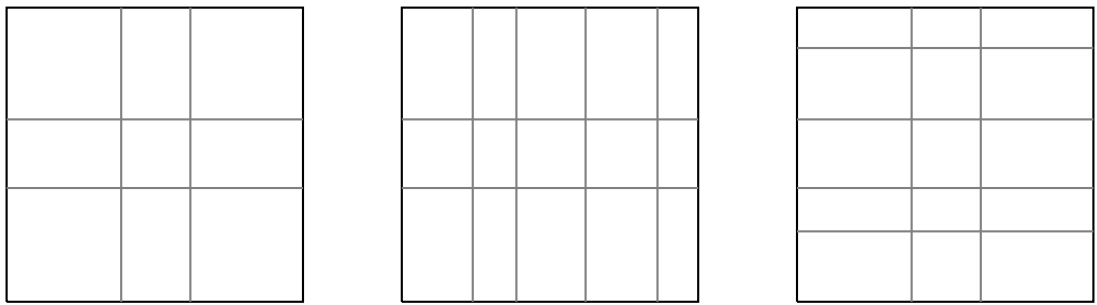
\begin{figure}[h!]
     \centering
    \begin{tikzpicture}[scale=1.3]

        \draw[thick, black] (0,0) -- (3,0) -- (3,3) -- (0,3) -- (0,0);
        \draw[thick, black] (4,0) -- (7,0) -- (7,3) -- (4,3) -- (4,0);
        \draw[thick, black] (8,0) -- (11,0) -- (11,3) -- (8,3) -- (8,0);

        \draw[thick, gray] (1.86,0) -- (1.86,3);
        \draw[thick, gray] (1.16,0) -- (1.16,3);
        \draw[thick, gray] (0,1.86) -- (3,1.86);
        \draw[thick, gray] (0,1.16) -- (3,1.16);
        
        \draw[thick, gray] (5.86,0) -- (5.86,3);
        \draw[thick, gray] (5.16,0) -- (5.16,3);
        \draw[thick, gray] (4,1.86) -- (7,1.86);
        \draw[thick, gray] (4,1.16) -- (7,1.16);
        \draw[thick, gray] (4.718,0) -- (4.718,3);
        \draw[thick, gray] (6.588,0) -- (6.588,3);

        \draw[thick, gray] (9.86,0) -- (9.86,3);
        \draw[thick, gray] (9.16,0) -- (9.16,3);
        \draw[thick, gray] (8,1.86) -- (11,1.86);
        \draw[thick, gray] (8,1.16) -- (11,1.16);
        \draw[thick, gray] (8,0.718) -- (11,0.718);
        \draw[thick, gray] (8,2.588) -- (11,2.588);

    \end{tikzpicture}
    \caption{The $(2,2)$, $(3,2)$ and $(2,3)$ partitions in base $\varphi$ from left to right. This should be used as an aid when reading Lemmas \ref{lem: point distribution} and \ref{lem: Extending points}.}
    \label{fig:22 32 23 partitions}
\end{figure}

\pagebreak
\newpage

\section{Numerical Discrepancy Results}\label{app:additional figures - discrepancy}

As further justification of the improvement in performance (with respect to the discrepancy) of digital nets and sequences constructed with irrational bases, we present several additional graphs and tables containing numerical studies of the star and $L_2-$discrepancy.

\begin{figure}[h!]
     \centering
     \includegraphics[scale=.5]{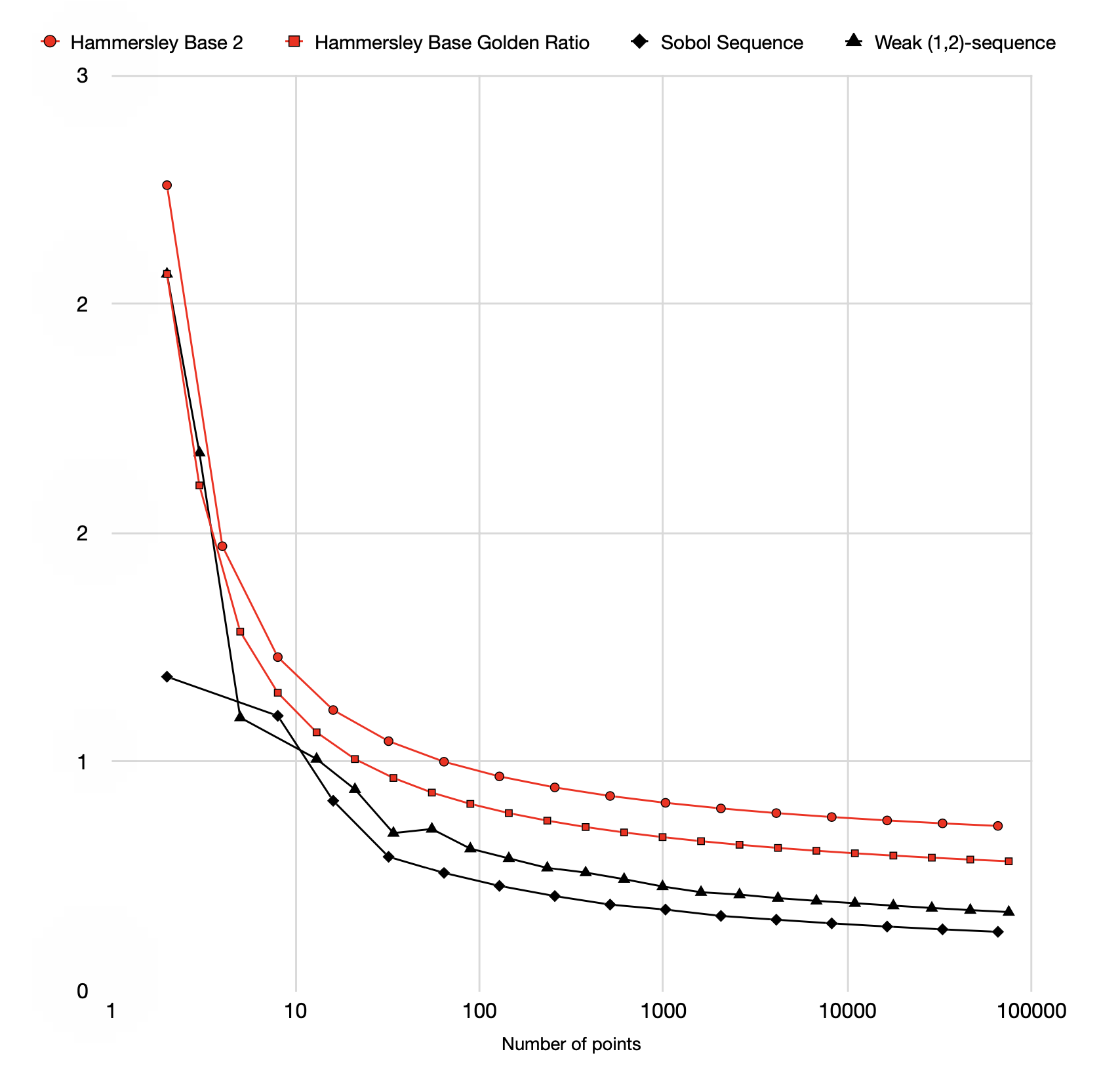}
     \caption{The $L_2-$discrepancy normalized by $\log N/N$ of the Hammersley point sets in base $2$ and $\varphi$, the weak $(1,2)-$sequence and the Sobol' sequence in two-dimensions.}\label{fig:L2 star disc}
\end{figure}

\newpage

To illustrate the phenomenon as described in Section \ref{sec:equidist in other bases} that the most well-distributed two-dimensional point set constructions are constructed from $\gamma$ derived when $q=1$ for any given $p$, we present the series of tables containing the star discrepancy values for a wide range of $N$. All values have been normalized by $\log N/N$.
\begin{table}[h!]
    \begin{minipage}{.5\linewidth}
      \centering
        \begin{tabular}{|| c | c ||}
            \hline
            $N$ & $D^*_N$ \\
            \hline \hline
            $3$ & $3.11$ \\
            $7$ & $2.60$ \\
            $17$ & $2.07$ \\
            $41$ & $1.75$ \\
            $99$ & $1.62$ \\
            $239$ & $1.44$ \\
            $577$ & $1.37$ \\
            $1393$ & $1.30$ \\
            $3363$ & $1.22$ \\
            $8119$ & $1.20$ \\
            $19601$ & $1.15$ \\
            \hline
            
        \end{tabular}
        \caption{$p=2,q=1$}\label{table:p2q1first}
    \end{minipage}%
    \begin{minipage}{.33\linewidth}
      \centering
        \begin{tabular}{|| c | c ||}
            \hline
            $N$ & $D^*_N$ \\
            \hline \hline
            $3$ & $3.35$ \\
            $8$ & $2.29$ \\
            $22$ & $1.92$ \\
            $60$ & $1.65$ \\
            $164$ & $1.61$ \\
            $448$ & $1.44$ \\
            $1224$ & $1.33$ \\
            $3344$ & $1.32$ \\
            $9136$ & $1.27$ \\
            $24960$ & $1.22$ \\
            $68192$ & $1.22$ \\
            \hline
  
        \end{tabular}
        \captionsetup{width=4cm}
        \caption{$p=2,q=2$}
    \end{minipage} 
\end{table}

 \begin{table}[h!]
    \begin{minipage}{.33\linewidth}
      \centering
        \begin{tabular}{|| c | c ||}
            
            \hline
            $N$ & $D^*_N$ \\
            \hline \hline
            $4$ & $2.71$ \\
            $13$ & $2.25$ \\
            $43$ & $1.83$ \\
            $142$ & $1.53$ \\
            $469$ & $1.47$ \\
            $1549$ & $1.36$ \\
            $5116$ & $1.31$ \\
            $16897$ & $1.30$ \\
            $55807$ & $1.26$ \\
            $184318$ & $1.20$ \\
            \hline
            
        \end{tabular}
        \captionsetup{width=4cm}
        \caption{$p=3,q=1$}
    \end{minipage}%
    \begin{minipage}{.33\linewidth}
      \centering
        \begin{tabular}{|| c | c ||}

            \hline
            $N$ & $D^*_N$ \\
            \hline \hline
            $4$ & $2.89$ \\
            $14$ & $2.20$ \\
            $50$ & $1.80$ \\
            $178$ & $1.55$ \\
            $634$ & $1.48$ \\
            $2258$ & $1.38$ \\
            $8042$ & $1.35$ \\
            $28642$ & $1.28$ \\
            $102010$ & $1.27$ \\
            \hline
            
        \end{tabular}
        \captionsetup{width=4cm}
        \caption{$p=3,q=2$}
    \end{minipage} 
    \begin{minipage}{.33\linewidth}
      \centering
        \begin{tabular}{|| c | c ||}
            \hline
            $N$ & $D^*_N$ \\
            \hline \hline
            $4$ & $3.13$ \\
            $15$ & $2.26$ \\
            $57$ & $1.94$ \\
            $216$ & $1.68$ \\
            $819$ & $1.62$ \\
            $3105$ & $1.53$ \\
            $11772$ & $1.49$ \\
            $44631$ & $1.44$ \\
            $169209$ & $1.41$ \\
            \hline
        \end{tabular}
        \captionsetup{width=4cm}
        \caption{$p=3,q=3$}
    \end{minipage} 
\end{table}

 \begin{table}[h!]
    \begin{minipage}{.24\linewidth}
      \centering
        \begin{tabular}{|| c | c ||}
            
             \hline
            $N$ & $D^*_N$ \\
            \hline \hline
            $5$ & $2.70$ \\
            $21$ & $2.13$ \\
            $89$ & $1.78$ \\
            $377$ & $1.56$ \\
            $1597$ & $1.47$ \\
            $6765$ & $1.40$ \\
            $28657$ & $1.35$ \\
            $121393$ & $1.30$ \\
            $514229$ & $1.27$ \\
            \hline
            
        \end{tabular}
        \captionsetup{width=4cm}
        \caption{$p=4,q=1$}
    \end{minipage}%
    \begin{minipage}{.24\linewidth}
      \centering
        \begin{tabular}{|| c | c ||}

            \hline
            $N$ & $D^*_N$ \\
            \hline \hline
            $5$ & $2.85$ \\
            $22$ & $2.11$ \\
            $98$ & $1.82$ \\
            $436$ & $1.55$ \\
            $1940$ & $1.52$ \\
            $8632$ & $1.40$ \\
            $38408$ & $1.38$ \\
            $170896$ & $1.32$ \\
            \hline
            
        \end{tabular}
        \captionsetup{width=4cm}
        \caption{$p=4,q=2$}
    \end{minipage} 
    \begin{minipage}{.24\linewidth}
      \centering
        \begin{tabular}{|| c | c ||}

            \hline
            $N$ & $D^*_N$ \\
            \hline \hline
            $5$ & $2.97$ \\
            $23$ & $2.03$ \\
            $107$ & $1.85$ \\
            $497$ & $1.65$ \\
            $2309$ & $1.57$ \\
            $10727$ & $1.50$ \\
            $49835$ & $1.47$ \\
            $231521$ & $1.42$ \\
            \hline

        \end{tabular}
        \captionsetup{width=4cm}
        \caption{$p=4,q=3$}
    \end{minipage} 
     \begin{minipage}{.24\linewidth}
      \centering
        \begin{tabular}{|| c | c ||}

            \hline
            $N$ & $D^*_N$ \\
            \hline \hline
            $5$ & $3.06$ \\
            $24$ & $2.19$ \\
            $116$ & $1.96$ \\
            $560$ & $1.77$ \\
            $2704$ & $1.73$ \\
            $13056$ & $1.63$ \\
            $63040$ & $1.61$ \\
            $304384$ & $1.56$ \\
            \hline
     
        \end{tabular}
        \captionsetup{width=4cm}
        \caption{$p=4,q=4$}\label{tab:p4q4last}
    \end{minipage} 
\end{table}

\end{document}